\documentclass[seceqn,secthm]{elsart}
\usepackage{latexsym,epsf,amssymb,multicol,makeidx}
\usepackage[dvips]{graphicx}
\usepackage{epsfig}

\usepackage{psfrag} %sustituye las etiquetas de los gr\'{a}ficos eps por etiquetas con formato igual al del texto

\usepackage{amsmath}
 \usepackage{capt-of}

\newcommand{\bq}{\begin{equation}}
\newcommand{\eq}{\end{equation}}
\newtheorem{theorem}{Theorem}
\newtheorem{lemma}{Lemma}

\newtheorem{definition}{Definition}
\newtheorem{remark}{Remark}

\newcommand{\norm}[1]{\left\Vert#1\right\Vert}

\begin{document}

\begin{frontmatter}


% use the corauthref command within \author for corresponding author
%footnotes;

\title{Error estimates for splitting methods  based on AMF-Runge-Kutta formulas for the time
integration of advection diffusion reaction PDEs.}

\author{
S. Gonzalez-Pinto,  D. Hernandez-Abreu and S. Perez-Rodriguez}

% $^{\mbox{\tiny {1}}}$

\thanks[thanks1]{This work has been supported by projects
MTM2010-21630-C02-02 and MTM2013-47318-C2-2}

\corauth[cor1]{Corresponding author: S. Gonzalez-Pinto
(spinto@ull.es)}

\baselineskip=0.9\normalbaselineskip \vspace{-3pt}

\maketitle

\address{
  { \footnotesize Departamento de An\'{a}lisis Matem\'{a}tico. Universidad de
La Laguna.\\
38071. La Laguna, Spain. \\
  email: spinto\symbol{'100}ull.es, dhabreu\symbol{'100}ull.es}}

\begin{abstract}

The convergence of a family of AMF-Runge-Kutta methods (in short
AMF-RK) for the time integration of evolutionary Partial
Differential Equations (PDEs) of Advection Diffusion Reaction type
semi-discretized in space is considered. The methods are based on
very few inexact Newton Iterations of Aproximate Matrix
Factorization splitting-type (AMF) applied to the Implicit
Runge-Kutta formulas, which allows very cheap and inexact
implementations of the underlying Runge-Kutta formula. Particular
AMF-RK methods based on Radau IIA formulas are considered. These
methods have given very competitive results when compared with
important formulas in the literature for multidimensional systems of
non-linear parabolic PDE problems. Uniform bounds for the global
time-space errors on semi-linear PDEs when simultaneously the time
step-size and the spatial grid resolution tend to zero are derived.
Numerical illustrations supporting the theory are presented.
\end{abstract}


\end{frontmatter}

%%%%%%%%%%%%%%%%%%%%%%%%%%%%%
%\input{sect1-PRE}

\section{Introduction}

We consider  numerical methods for the time integration of a family
of Initial Value Problems in ODEs
\begin{equation}\label{ode}
y_h'(t) = f_h(t, y_h(t)),\;\;\;y_h(0) = {u}^*_{0,h}, \;\;\; 0 \le t
\le t^*, \quad y_h, f_h \in \mathbb{R}^{m(h)}, \quad h\rightarrow
0^+,
\end{equation}
coming from the spatial semi-discretization of an $l-$dimensional
 Advection Diffusion Reaction  problem  in time
dependent Partial Differential Equations (PDEs), with prescribed
Boundary Conditions and an Initial Condition.  Here $h$ denotes  a
small positive parameter associated with the spatial resolution and usually $l=2,3,\ldots$ .

The typical PDE problem with Dirichlet
boundary conditions is given by ($\Omega$ is
a bounded open connected  region in $\mathbb{R}^l$, $\partial \Omega$  its boundary and $\nabla$ is the
gradient operator)
\begin{equation}\label{pde} \begin{array}{c}
u_t(x,t) = - \nabla \cdot (a(x,t) u(x,t)) + \nabla \cdot
(\bar{d}(x,t)\cdot \nabla u(x,t)) + r(u,x,t),\\[0.5pc]  x \in \Omega , \; t\in [0,t^*]; \; a(x,t)=(a_j(x,t))_{j=1}^l \in \mathbb{R}^l, \;\bar{d}(x,t)=(\bar{d}_j(x,t))_{j=1}^l \in \mathbb{R}^l,\\[0.3pc]
u(x,t) =g_1(x,t), \:(x,t)\in  \partial \Omega\times[0,t^*];
\qquad u(x,0)=g_2(x), \; x \in \Omega,\end{array}
\end{equation}
which is  assumed to have some diffusion ($\bar{d}_j(x,t)\ge
d_0>0,\;j=1,\ldots,l$), namely that it is not of pure hyperbolic
type, and  it is also assumed that some adequate spatial
discretization based on Finite Differences or Finite Volume is
applied to obtain the system (\ref{ode}).   Some stiffnes in the
reaction part $r(u,x,t)$ is also allowed. The treatment of Systems
of PDEs do not involve more difficulty for our analysis  but for
simplicity of presentation we prefer to confine ourselves to the
case of one PDE.

We denote by $u_h(t)$ the solution of the PDE problem confined to
the spatial grid (or well to the $h$-space related). It will be
tacitly assumed that the PDE problem admits a smooth solution
$u(x,t)$ in the sense that  continuous partial derivatives in all variables  up to some order $p$
exist and  are continuous and uniform bounded on $\Omega\times[0,t^*]$ and that $u(x,t)$ is continuous
on  $\bar{\Omega}\times[0,t^*]$ ($\bar{\Omega}=\Omega \bigcup \partial \Omega$).
It is also assumed that the spatial discretization errors \begin{equation}\label{spatialerrors}
\sigma_h(t):= u_h'(t)-f_h(t,u_h(t)),
\end{equation}
 satisfy in the norm considered,
\begin{equation}\label{norm}%
\Vert \sigma_h(t) \Vert \le C \:h^r , \quad (C\ge 0,\;r>0), \quad 0
\le t \le t^*, \quad h\rightarrow 0.
\end{equation}
In general $C$, $C'$ or $C^*$ will refer to some constants that
maybe different at each occurrence but that all of them are
independent of $h\rightarrow 0$ and from the time-stepsize
$\tau\rightarrow 0$. The vector norm used is arbitrary as long as it
is defined for vectors of any dimension. For square matrices the
norm used is the induced operator norm, $\Vert A\Vert = \sup_{v\ne
0} \Vert A v \Vert / \Vert v \Vert. $

In spite of  most of our results apply  in general, we will provide
specific results for  weighted Euclidean norms of type
$$\Vert (v_j)_{j=1}^N \Vert=N^{-1/2} \Vert (v_j)_{j=1}^N \Vert_2. $$
It should be noted that in this case we have for any square matrix
$A$ that,
$$  \Vert A \Vert=  \Vert A \Vert_2, \quad \forall \:A \in \mathbb{R}^{N,N},\;N=1,2,3,\ldots .$$
We  assume some natural
splitting for $f_h$ (directional or other),
\begin{equation}\label{splitting}
f_h(t,y)=\sum_{j=1}^d f_{j,h}(t,y),
\end{equation}
which provides some natural splitting for the Jacobian matrix at the current point $(t_n,y_n)$,
\begin{equation}\label{split}
J_h=\sum_{j=1}^d J_{j,h}, \quad J_h :=
\displaystyle{\frac{\partial f_h(t_n,y_n)}{\partial
y}}, \quad J_{j,h} :=\displaystyle{\frac{\partial f_{j,h}(t_n,y_n)}{\partial y}}.
\end{equation}

This  goal of the paper is to analyze  the convergence order  of
 the Method of Lines (MoL) approach  for time-dependent PDEs of Advection Reaction Diffusion
 PDEs, with the main focuss on the time integration of the large ODE systems
resulting of the spatial PDE-semidiscretization, where some
stiffness is assumed (parabolic dominant problems with stiff
reaction terms) and the time integrators are based on very few
iterations of  splitting type (Approximate Matrix Factorization and
Newton-type schemes) applied to highly stable Implicit Runge-Kutta
methods. It should be remarked that the underlying Implicit
Runge-Kutta method is never solved up to convergence, hence the
convergence study does not follows from the results collected in
classical references about finite difference methods such as
\cite{Rich-Morton67,Bur-Hun-Ver86,Marchuk90,Thomee90,Trefethen92,HV}.
The kind of approach to be considered here has interest since it is
easily applicable to general systems of PDEs as we will see later on
and it is
 reasonably cheap for non-linear problems in general
 (although we give convergence results for semilinear problems only)  when some  splitting of the function $f_h$ and
 its Jacobian is available and the split terms can be  handled efficiently. In particular a method based on
 three AMF-iterations of the two-stage Radau IIA method
 \cite{Axelsson69} has shown to be competitive \cite{apnum-sevsole10} when compared with some
 standard PDE-solvers  such as VODPK \cite{Brown-Byrne-Hindmarsh-SISC89,Brown}
 in some interesting non-linear diffusion reaction problems
widely  considered in the literature. We also present two new
methods based on the 2-stage Radau IIA, by performing just  one or
two iterations of splitting type, respectively.  The method based on
two iterations is one of the very few one-step methods of splitting
type we have seen in the literature that has order three in
PDE-sense for the time integration.

The rest of the paper is organized as follows. In section 2 we
introduce the {\sf AMF$_q$-RK} methods, and  special attention is
paid to some methods based on Radau IIA formulas. In section 3, the
convergence for semilinear PDEs is studied in detail. The local and
global errors are studied for  the {\sf AMF$_q$-RK} splitting
methods based on some general Runge-Kutta methods.  Section 4 is
devoted to some applications of the convergence results to 2D and
3D-parabolic PDEs.

 {\rm Henceforth, for simplicity in the
notations,  we  omit in many cases  the $h$-dependence
 of some vectors such as  $f_h,\;f_{j,h}$ and of some matrices such as $J_h$
 and $J_{j,h}$ ($j=1,\ldots,d$). It should be clear from the context which ones are $h$-dependent.
 Besides, we will refer to the identity matrix as $I$ when its dimension is clear from the context. }

%%%%%%%%%%%%%%%%%%%%%%%%%%%%%%%%%%%%%%
%\input{sect2-PRE}
\section{AMF-IRK methods}

For the integration of  the ODEs (\ref{ode}), we consider as a first
step an implicit s-stage Runge-Kutta method with a nonsingular
coefficient matrix  $A=(a_{ij})_{i,j=1}^s$ and a weight vector
$b=(b_j)_{j=1}^s$. The method is  given by the  compact formulation
(below $\otimes$ denotes the Kronecker product of matrices $A\otimes
B=(a_{ij}B), \;A=(a_{ij}), \;B=(b_{ij})$)
\begin{equation}\label{IRK} \begin{array}{c}
Y_n=e\otimes y_n +\tau (A\otimes I_m)F(Y_n), \\
y_{n+1}= \varpi y_n + (\ss^T\otimes I_m)Y_n, \\ c\equiv (c_j)_{j=1}^s:=A e, \quad e=(1,\ldots,1)^T \in  \mathbb{R}^{s}, \quad \ss^T:=b^TA^{-1},\quad \varpi=1-\ss^Te, \\
Y_n=(Y_{n,j})_{j=1}^s \in \mathbb{R}^{ms}, \qquad  F(Y_n)=(f(t_n+\tau c_j,Y_{n,j}))_{j=1}^s \in \mathbb{R}^{ms}.
\end{array}
\end{equation}
It should be noted that we have replaced the usual formulation at the stepping point $y_{n+1}= y_n + \tau (b^T\otimes I_m) F(Y_n)$ by the equivalent in (\ref{IRK}), which has some computational advantages for stiff problems when the algebraic system for the stages is not exactly solved.

A typical Quasi-Newton iteration to solve the stage equations above
is given by (below, $J=\partial{f}/\partial y\:(t_n,y_n)$ is the
exact Jacobian at the step-point $(t_n,y_n)$),
\begin{equation}\label{newt}
[I_{ms}- A\otimes \tau J]\Delta^\nu=D_n^{\nu-1},\quad
Y_n^\nu=Y_n^{\nu-1}+\Delta^\nu, \quad \nu=1,2,\dots,
\end{equation}
where
\begin{equation}\label{residual}
D_n^{\nu-1}\equiv D(t_n,\tau,y_n,Y_n^{\nu-1}):=
e\otimes y_n - Y_n^{\nu-1} + \tau ( A\otimes I_m)
F(Y_n^{\nu-1}).
\end{equation}
A cheaper iteration of Newton-type   when the matrix $A$ has a
multipoint spectrum  has been considered in
\cite{APNUM95-seve,PGS09} (denoted as Single-Newton iteration)
\begin{equation}\label{s-newt}
[I_{ms}- T_\nu\otimes \tau J]\Delta^\nu=D_n^{\nu-1},\quad
Y_n^\nu=Y_n^{\nu-1}+\Delta^\nu, \quad \nu=1,2,\dots,q
\end{equation}
where
\begin{equation}\label{T} \begin{array}{c}
T_\nu=\gamma S_\nu (I-L_\nu)^{-1} S_\nu^{-1}, \quad \gamma>0,\\
\quad S_\nu \in \mathbb{R}^{s,s} \;\mbox{\rm are  regular matrices and }
\\ L_\nu \in \mathbb{R}^{s,s} \;   \mbox{\rm are strictly lower triangular matrices.}
\end{array}
\end{equation}
After some simple manipulations, by using standard properties of the Kronecker product,
this iteration can be rewritten in the equivalent form,
\begin{equation}\label{s-newt1}
\begin{array}{rl}
[I_s\otimes(I_{m}-\gamma \tau J)] E^\nu&=((I_s-L_\nu)S_\nu^{-1} \otimes I_m)D_n^{\nu-1}+ (L_\nu\otimes I_m)E^\nu, \\
Y_n^\nu &=Y_n^{\nu-1}+(S_\nu\otimes I_m)E^\nu, \qquad \nu=1,2,\dots,q.
\end{array}
\end{equation}
To reduce the algebra cost, we use the Approximate Matrix
Factorization \cite{How-Som-JCAM2001} in short AMF, with $J\equiv
J_h$ and $J_j\equiv J_{j,h}$ given in (\ref{split}),
\begin{equation}\label{product}
\Pi_d := \prod_{j=1}^d (I_m-\gamma \tau J_j)= (I_{m}-\gamma \tau J) +
\mathcal{O}(\tau^2),
\end{equation}
and replace in (\ref{s-newt1})  $(I_{m}-\gamma \tau J)$ by $\Pi_d$,
which yields the {\sf AMF$_q$-RK} method based on the underlying
Runge-Kutta method
\begin{equation}\label{AMF-RK}
\begin{array}{rl}
[I_s\otimes\Pi_d] E^\nu&=((I_s-L_\nu)S_\nu^{-1} \otimes I_m)D_n^{\nu-1}+ (L_\nu\otimes I_m)E^\nu, \\
Y_n^\nu &=Y_n^{\nu-1}+(S_\nu\otimes I_m)E^\nu, \qquad \nu=1,2,\dots,q \\
Y_n^0&=e\otimes y_n \qquad \mbox{\rm (Predictor)} \\ y_{n+1}&= \varpi y_n + (\ss^T\otimes I_m)Y^q_n
\qquad \mbox{\rm (Corrector).}
\end{array}
\end{equation}
Our starting point for the convergence analysis in the next section
takes into account that the {\sf AMF$_q$-RK} method can be rewritten
in the equivalent form \cite{sevedom-AMFestab}

\begin{equation}\label{AMF-RK1} \begin{array}{c}
 [I\otimes I-
T_\nu\otimes \tau P](Y_n^\nu-Y_n^{\nu-1})=
D(t_n,\tau,y_n,Y_n^{\nu-1}), \;\; 1\leq \nu\leq q \\
Y_n^0=e\otimes y_n, \qquad  y_{n+1}= \varpi y_n + (\ss^T\otimes I_m)Y^q_n,
\end{array}
\end{equation}
where the matrix $P$ plays a primary role
\begin{equation}\label{matrixP}
\begin{array}{lll} P&:=&(\gamma
\tau)^{-1}(I-\Pi_d)\\&=& J+(-\gamma
\tau)\displaystyle{\sum_{j<k}} J_jJ_k +(-\gamma
\tau)^2 \displaystyle{\sum_{j<k<l}} J_jJ_kJ_l
+\ldots+(-\gamma \tau)^{d-1} J_1J_2\cdots
J_d.\end{array}
\end{equation}
\subsection{AMF$_q$-RK methods based on the 2 stage Radau IIA
formula}\label{sec-2.1}

We are going to deserve special attention to  AMF$_q$-RK methods
 based on the 2 stage Radau IIA formula \cite{Axelsson69}. This  formula  has
coefficient Butcher tableau given by
$$\begin{array}{c|c} c&A\\ \hline \\[-2pc] & b^T\end{array} \quad \equiv
\quad \begin{array}{c|cc} 1/3& 5/12 & \;-1/12\\   1 &  3/4&\; 1/4
\\[0.3pc] \hline \\[-2pc]
  &  3/4&\; 1/4\end{array}$$
  This is a collocation method (stage order is two) possessing good stability properties, such as
$L$-stability (i.e. $A$-stability plus $R(\infty)=0$, with $R(z)$
being the linear stability function of the method), and has order of
convergence three (in ODE sense), not only on non-stiff problems but
also in many kinds of stiff problems \cite{Bur-Hun-Ver86}. These
properties for the underlying Runge-Kutta method are convenient,
since the family of ODEs (\ref{ode}) involves stiffness in most of
cases, due to the diffusion terms and possibly to the reaction part,
and it is expected that the methods to be built on inherit part of
the good properties of the original Runge-Kutta method.

The next three {\sf AMF$_q$-Rad} methods have coefficient matrices
($L_\nu$, $S_\nu$ and $T_\nu$) and eigenvalue $\gamma$ of the form
\begin{equation}\label{T}
T_\nu=\gamma S_\nu (I_2-L_\nu)^{-1} S_\nu^{-1}, \;
S_\nu=\left(\begin{array}{cc} 1 & s_\nu
\\ 0&1\end{array}\right),\; L_\nu=\left(\begin{array}{cc} 0 & 0 \\
l_\nu  &0\end{array}\right),\; \gamma=\sqrt{\det(A)}=1/\sqrt{6}.
\end{equation}
{\sf AMF$_1$-Rad} was derived in \cite{sevedom-AMFestab} by looking
for  good stability properties and order two (ODE sense). In
particular the method is A($\pi/2$)-stable for a $2$-splitting (see
in Definition \ref{def-estability} below, the concept of stability
for a $d$-splitting), A($0$)-stable for any $d$-splitting and has
stability wedges close to $\theta_d=\pi/(2(d-1))$ for $d=3,4$. The
method is based on one iteration ($q=1$) and was required to fulfil
$(A-T_1)c=0$ and it has coefficients given by
\begin{equation}\label{met1}
s_1= -\frac{3+2 \sqrt{6}}{9}, \quad l_1=\frac{3}{4}(-12+5\sqrt{6}).
\end{equation}

{\sf AMF$_2$-Rad} was derived in \cite{sevedom-AMFestab} by looking
for  good stability properties and order three (ODE sense). The
method is A($\pi/2$)-stable for a $2$-splitting, A($0$)-stable for
any $d$-splitting and  A($\pi/6$)-stable  for $d=3,4$. The method is
based on two iterations ($q=2$) and their matrices $T_1$ and $T_2$
were required to satisfy $(A-T_1)c=0$ and $e_2^T
T_2^{-1}(A-T_2)=0^T, \;e_2^T=(0,1)$, respectively.  Its coefficients
are uniquely given by
\begin{equation}\label{met2} \begin{array}{c}
\displaystyle{s_1= -\frac{3+2 \sqrt{6}}{9}, \quad
l_1=\frac{3}{4}(-12+5\sqrt{6})}\\[0.3pc]
\displaystyle{s_2= \frac{5-2\sqrt{6}}{9}, \quad
l_2=\frac{3\sqrt{6}}{4}}.\end{array}
\end{equation}

{\sf AMF$_3$-Rad} was derived in \cite{PGS09,apnum-sevsole10} by
looking for good stability properties and order three (ODE sense).
The method is A($\pi/2$)-stable for a $2$-splitting, A($0$)-stable
for any $d$-splitting and  close to A($\theta_d$)-stable  for
$d=3,4$ with $\theta_d=\pi/(2(d-1))$. The method is based on three
iterations ($q=3$) and their matrices $T=T_1=T_2=T_3$  were required
to satisfy $e_2^T T^{-1}(A-T)=0^T$. Its coefficients are uniquely
given by
\begin{equation}\label{met3} \begin{array}{c}
\displaystyle{s_1=s_2=s_3= \frac{5-2\sqrt{6}}{9}, \quad
l_1=l_2=l_3=\frac{3\sqrt{6}}{4}}.\end{array}
\end{equation}
 In \cite{apnum-sevsole10},  a variable-stepsize integrator based
on the {\sf AMF$_3$-Rad} method  was successfully tested on several
interesting $2D$ and $3D$ advection diffusion reaction PDEs by
exhibiting good performances in comparison with state-of-the-art
codes like {\sf VODPK} \cite{Brown-Byrne-Hindmarsh-SISC89,Brown} and
{\sf RKC} \cite{RKC,VSH-JcompPhys2004} and its implicit-explicit
counterpart, {\sf IRKC} \cite{IRKC,IMEXRKC}. The other two methods,
{\sf AMF$_q$-Rad} ($q=1,2$), were introduced later
\cite{sevedom-AMFestab} after carefully analyzing the PDE errors  on
semilinear problems and with the purpose of reducing the number of
iterations w.r.t. {\sf AMF$_3$-Rad}.
%%%%%%%%%%%%%%%%%%%%%%%%%%%%%%%%%%%%%%
%\input{sect2-PRE}
\section{Convergence for semilinear problems}
For our convergence analysis we consider   {\sf AMF$_q$-RK} methods
 applied to the ODE problems coming from the spatial
discretizations of semilinear PDE problems of type (\ref{pde}) where
the advection and diffusion vectors $a(x,t)$ and $\bar{d}(x,t)$ are
both constant and the reaction part has the form
\begin{equation}\label{reaction} r(u,x,t)=\kappa\: u + g(x,t), \quad \kappa \;
\mbox{\rm being a constant},\quad x\in\Omega \subseteq \mathbb{R}^l.
\end{equation} In this way, the ODE systems have the form
\begin{equation}\label{lin-system}
\begin{array}{c}
y_h'(t)=f_h(t,y_h):=J_h y_h(t)+g_h(t),\quad
y_h(0)=u^*_{0,h},\quad h\rightarrow 0^+,\\
J_h=\sum_{j=1}^d
J_{j,h}, \qquad t\in [0,t^*].
\end{array}
\end{equation}
Here, the  exact solution of the PDE  confined to the spatial grid
$u_h(t)=u(x,t)$, is assumed to satisfy (\ref{spatialerrors}) and
(\ref{norm}). Thus, we focus on the global errors of the MoL
approach, where the spatial discretization is carried out first by
using finite differences (or finite volumes) and then the time
discretization is performed by  using {\sf AMF$_q$-RK} methods. It
is important to remark that we will not pursue the details of the
spatial semidiscretizations but rather it is assumed that the
spatial semidiscretizations  are stable and provides spatial
discretization errors satisfying  (\ref{norm}). We shall provide
uniform bounds for the global errors of the MoL approach ($y_h(t)$
henceforth denotes the numerical solution of the MoL approach) in
the sense
\begin{equation}\label{global-errors}
\epsilon_{n,h}:=u_h(t_n)-{y}_h(t_n)=
\mathcal{O}(\tau)^{p_1}
+\mathcal{O}(h^{\alpha}\tau^{p_2}) ,\quad h\rightarrow 0^+,\tau
\rightarrow 0^+,
\end{equation}
which is meant that there exist  constants $C_1,\:C_2,\:p_1,\:p_2,\:\alpha$
(all of them independent on $h$ and $\tau$)   so that in the norm considered,
$$\Vert \epsilon_{n,h}\Vert \le C_1 \tau^{p_1} + C_2 h^{\alpha} \tau^{p_2},\quad h \rightarrow 0^+,\tau
\rightarrow 0^+ \quad \mbox{\rm holds. }$$ In our convergence
analysis we need that all the  matrices $J_{j,h}$ pairwise conmute
and that they can be brought to the following decomposition (it has
some resemblance with the Jordan's decomposition, but it is a little
more general)
\begin{equation}\label{jordan}
\begin{array}{c}
J_{j,h}= \Theta_{h} \Lambda_{j,h}
\Theta_{h}^{-1},  \quad \mbox{\rm Cond}
(\Theta_{h}):=\Vert \Theta_{h}\Vert\cdot \Vert
\Theta_{h}^{-1}\Vert \le C,
\;h\rightarrow 0^+, \;1\leq j\leq d,\\
\Lambda_{j,h}=\mbox{\rm
BlockDiag}(\Lambda^{(1)}_{j,h},\Lambda^{(2)}_{j,h},\ldots,\Lambda^{(\vartheta_h)}_{j,h}),
\quad  \Lambda^{(l)}_{j,h}=\lambda^{(l)}_{j,h} I + E^{(l)}_{h} , \quad \mbox{\rm Re }  \lambda^{(l)}_{j,h}  \le 0,\\
  \mbox{\rm dim}(E^{(l)}_{h})\le N, \quad \Vert E^{(l)}_{h} \Vert_\infty \le C',
  \; l=1,2,\ldots,\vartheta_h  \quad (h\rightarrow 0^+).\\
E^{(l)}_{h} \quad \mbox{\rm are all of them strictly lower triangular matrices. }
\end{array}
\end{equation}

Another important approach for the convergence analysis of the MoL
method (mainly concerned with the time integration) is based on the
pseudo-spectra analysis of the matrix $J_{h}$ \cite{Trefethen92} and
the related matrices $J_{j,h}$. That analysis is of more general
scope but it is much more difficult to make and as we will see
below, our analysis is enough for some interesting kind of
semilinear problems and it is expected that the results extend to
most of the non-linear problems of parabolic dominant type.

Next, we consider a standard 3D semilinear-PDEs problem where the
assumptions in (\ref{jordan}) are fulfilled.

\subsection{An example}
Consider the semilinear PDE-problem (\ref{pde}) with $x\in
\Omega=(0,1)^3$, with constant vectors, $a(x,t)=(a_j)_{j=1}^3, \;
\bar{d}(x,t)=(\bar{d}_j)_{j=1}^3, \;\bar{d}_j>0 \;(j=1,2,3)$ and
$r(x,u,t)$ as in (\ref{reaction}). Consider the spatial
semidiscretization by using second order central differences and
spatial resolution $h=1/(N+1)$. This yields a semilinear ODE systems
of dimension $m=N^3$ of the form (\ref{lin-system}) for $d=3$. The
matrices $J_{j,h}$ are given by
\begin{equation}\label{J-example} \begin{array}{c}
J_{1,h}=I_N\otimes I_N\otimes \mathcal{T}_1,\quad J_{2,h}=I_N\otimes
\mathcal{T}_2\otimes I_N,\quad J_{3,h}= \mathcal{T}_3\otimes I_N
\otimes I_N\\ \mathcal{T}_l =\mbox{\rm
Tridiag}(\alpha_l,\delta_l,\beta_l)\in \mathbb{R}^{N,N},\quad l=1,2,3,\\
\alpha_l=h^{-2}( \bar{d}_l- 2^{-1}h \: a_l ), \quad \beta_l=h^{-2}(
\bar{d}_l+ 2^{-1}h \:a_l ), \quad \delta_l=h^{-2}( -2\bar{d}_l+h^2
\kappa),
\end{array}
\end{equation}
and the vector $g_h(t)$ includes the reaction part $g(x,t)$ plus the
boundary conditions. It is straightforward to see that the $J_{l,h}$
pairwise commute. Moreover, by assuming Cell-P\'{e}clet numbers \cite[p.
67, formula (3.42) ]{HV}
$$ h|a_l|/\bar{d}_l< 2 ,\quad l=1,2,3,$$  from \cite[section
2]{Pasquini13} it follows that their spectral decomposition has the
form
\begin{equation}\label{Tj} \begin{array}{c}
 \mathcal{T}_l =\mbox{\rm
Tridiag}(\alpha_l,\delta_l,\beta_l) =V_l \Lambda_lV_l^{-1},\quad V_l= D_l U, \quad l=1,2,3,\\[0.2pc]
\Lambda_l=\mbox{\rm Diag}\displaystyle{\left(
\lambda_{l,k}\right)_{k=1}^N,\quad \lambda_{l,k}= \delta_l +
2\sqrt{\alpha_l\beta_l} \cos{\frac{k\pi}{N+1}},}
\\[0.2pc] U=(\frac{2}{N+1})^{1/2}\displaystyle{\left(\sin{\frac{kj\pi}{N+1}}\right)_{k=1,N
\atop j=1,N}}
\; \mbox{\rm is an orthogonal matrix and } \\[0.2pc]
D_l= (\frac{N+1}{2})^{1/2}\mbox{\rm
Diag}\displaystyle{\left((\alpha_l/\beta_l)^{k/2}\right)_{k=1}^N.}
\end{array}
\end{equation}
From here we conclude that all the matrices can be brought to the
spectral decomposition in   (\ref{jordan}) having negative
eigenvalues and with matrix $\Theta_{h}=V_3\otimes V_2\otimes V_1. $
Observe that
$$\begin{array}{rcl} \Vert \Theta_{h}\Vert_2\Vert
\Theta_{h}^{-1}\Vert_2&=& \prod_{l=1}^3 \Vert V_l \Vert_2 \Vert
V_l^{-1} \Vert_2= \prod_{l=1}^3 \Vert D_l \Vert_2 \Vert D_l^{-1}
\Vert_2\\ &=& \displaystyle{ \prod_{l=1}^3 \left(\frac{2\bar{d}_l
+h|a_l|}{2\bar{d}_l -h|a_l|}\right)^{N/2}\le  \prod_{l=1}^3
\left(\frac{2\bar{d}_l +h|a_l|}{2\bar{d}_l -h|a_l|}\right)^{1/(2h)}}
\\ & \simeq & \displaystyle{\exp\left(\sum_{l=1}^3 \frac{|a_l|}{2\bar{d}_l}\right)}
\quad \mbox{\rm as } h\rightarrow 0.
\end{array}
$$

\subsection{Analysis of the Truncation Errors}
%For facility of presentation we will omit in some places the
%dependence on $h$.
The {\sf AMF$_q$-RK} method applied on problem (\ref{ode}) can be
expressed in the simple one-step format
$y_{n+1}=\phi_f(t_n,y_n,\tau),\;n\geq 0$. Thus, the time-space
global errors $\epsilon_{n}=u_h(t_n)-y_n$ satisfy
$$
\begin{array}{lll}\epsilon_{n+1}&:=&
u_h(t_{n+1})-\phi_f(t_n,y_n,\tau)\\
&=& (u_h(t_{n+1})-\phi_f(t_n,u_h(t_n),\tau)) +
(\phi_f(t_n,u_h(t_n),\tau)-\phi_f(t_n,y_n,\tau))\\
&=&l(t_n,\tau,h)+[\partial \phi_f/\partial y]_n (u_h(t_n)-y_n),
\end{array}
$$ where
%\begin{equation}\label{sev2-1}
$$
[\partial \phi_f/\partial y]_n=\int_0^1 \frac{\partial
\phi_f}{\partial
y}(t_n,u_h(t_n)+(\theta-1)\epsilon_{n},\tau)d\theta,
$$
%\end{equation}
and the {\it time-space local  errors} are defined by
\begin{equation}\label{local-errors}
l_n\equiv l(t_n,\tau,h) :=u_h(t_{n+1})-\phi_f(t_n,u_h(t_n),\tau).
\end{equation}
Then, we have for the {\it time-space global errors} $\epsilon_{n}$
the recurrence
\begin{equation}\label{global-error}
\epsilon_{n+1}= [\partial \phi_f/\partial y]_n \cdot \epsilon_{n}+
l_n,\quad n=0,1,2,\ldots,t^*/\tau-1.
\end{equation}
In order to get a better understanding of the latter recurrence, we
next introduce the following matrix operators ($P$ is defined in
(\ref{matrixP}))
\begin{equation}\label{Mnu-Qnu}
Q_\nu=(I\otimes I- T_\nu\otimes \tau P)^{-1}, \quad M_\nu=Q_\nu
(A\otimes \tau J - T_\nu \otimes \tau P),\quad   \nu\geq 1, \;Q_0=I.
\end{equation}

\begin{lemma}\label{GE-recursion}
The time-space global errors provided by the {\sf AMF$_q$-RK} method
 when applied to the problem (\ref{lin-system})
satisfy the recurrence
\begin{equation}\label{global-errors1}
\epsilon_{n+1}=  R_q(\tau J,\tau P)\cdot \epsilon_{n}+
l_n,\quad n=0,1,2,\ldots,t^*/\tau-1,
\end{equation}
where $l_n$ stands for the time-space local error defined in
(\ref{local-errors}) and
\begin{equation}\label{estab}\begin{array}{l}
R_q(\tau J,\tau P)= \varpi I + \displaystyle{(\ss^T\otimes
I)\left( Q_q + \sum_{j=q}^1 (\prod_{i=q}^{j}
M_i)Q_{j-1}\right) (e\otimes I)},
\end{array} \end{equation}
with  $Q_\nu, M_\nu$ ($\nu\geq 1$) given by
(\ref{Mnu-Qnu}). Moreover, the function $R_q(\tau J,\tau P)$ fulfils
\begin{equation}\label{estab-1} R_q(\tau J,\tau P)-I=  (\ss^T\otimes I)\left(Q_q + \sum_{j=q}^1
(\prod_{i=q}^{j} M_i)Q_{j-1}-\prod_{i=q}^{1} M_i\right) (c\otimes \tau J).
\end{equation}
\end{lemma}

\begin{remark} {\rm It must be observed that commutativity does not hold in
general, thus $\prod_{j=q}^1 M_j\equiv M_q M_{q-1}\cdots M_1.$
 On the other hand, $R_q(\cdot)$ can be
seen as the linear stability function of the method. The identity
(\ref{estab-1}) for the function $R_q(\cdot)-I$ will play a major
role in a favourable propagation of the local errors  in a similar
way as indicated in Lemma 2.3 in \cite[p.162]{HV}.}
\end{remark}

\noindent {\bf Proof of Lemma \ref{GE-recursion}.} Our first step is
to analyze the operator $[\partial \phi_f/\partial y]_n$ for the
semilinear problem (\ref{lin-system}). Taking into account that the
method is defined by (\ref{AMF-RK1}), then we are led to compute
$\displaystyle{\frac{\partial y_{n+1}}{\partial y_n}}$ with $
y_{n+1}= \varpi y_n + (\ss^T\otimes) Y_n^q$. At this end, by taking
derivatives with regard to $y_n$ in the iteration (\ref{AMF-RK1}),
it holds that
$$
\begin{array}{lll}
(I\otimes I- T_\nu\otimes \tau P)\left(\dfrac{\partial
Y_n^\nu}{\partial y_n}-\dfrac{\partial Y_n^{\nu-1}}{\partial
y_n}\right)&=& \dfrac{\partial
D(t_n,\tau,y_n,Y_n^{\nu-1})}{\partial y_n}\\
&=& e\otimes I +(-I\otimes I + A\otimes \tau J) \dfrac{\partial
Y_n^{\nu-1}}{\partial y_n}.
\end{array}
$$
From here, after some simple manipulations it follows that,
\begin{equation}\label{deriv}\begin{array}{l}
\displaystyle{\frac{\partial Y_n^\nu}{\partial y_n}}=
\displaystyle{Q_\nu (e\otimes I) + M_\nu \frac{\partial
Y_n^{\nu-1}}{\partial y_n}}, \quad (\nu=1,2,\ldots,q), \quad
\displaystyle{\frac{\partial Y_n^0}{\partial y_n}=e\otimes I.}
\end{array} \end{equation}
From an inductive argument, it is not difficult to see that
\begin{equation}\label{sol1}\begin{array}{l}
\displaystyle{\frac{\partial Y_n^q}{\partial y_n}}=
\displaystyle{\left( Q_q + \sum_{j=q}^1
(\prod_{i=q}^{j} M_i)Q_{j-1}\right) (e\otimes I).}
\end{array} \end{equation}
Then, by denoting $R_q(\tau J,\tau P):=\displaystyle{\frac{\partial
y_{n+1}}{\partial y_n}}$ it follows $$R_q(\tau J,\tau P)= \varpi I +
(\ss^T\otimes)\displaystyle{\frac{\partial Y_n^q}{\partial y_n}},$$
and we deduce both (\ref{estab}) and (\ref{global-errors1}) from
(\ref{global-error}) and (\ref{sol1}).

In order to prove (\ref{estab-1}),  we first take into account that
$R_q(\cdot)-I= (\ss^T\otimes I)Z_n^q$, where
 $Z_n^\nu= \partial Y_n^\nu/\partial y_n-e\otimes I$. Then, from the
recurrence (\ref{deriv}), it follows after some simple calculations
that $Z_n^\nu = M_\nu Z_n^{\nu-1} + Q_\nu (c\otimes \tau J)$,
$(\nu=1,2,\ldots,q)$, with $Z_n^0=0$. From here, we deduce
$Z_n^q=\displaystyle{\left( Q_q + \sum_{j=q}^1 (\prod_{i=q}^{j}
M_i)Q_{j-1}-\prod_{i=q}^{1}
M_i\right) (c\otimes \tau J)},$ and this directly gives
(\ref{estab-1}).\hfill $\Box$

\begin{remark}\label{sev-remark-0} {\rm  For a given rational function of two complex variables
\begin{equation}\label{sev-res0} \displaystyle{\zeta(z,w)=\frac{\sum_{i,j=0}^{m_1} \alpha_{ij}
z^iw^j}{\sum_{i,j=0}^{m_2} \beta_{ij}z^iw^j}\equiv  \left(
{\sum_{i,j=0}^{m_1} \alpha_{ij} z^iw^j}\right)
\left({\sum_{i,j=0}^{m_2}
\beta_{ij}z^iw^j}\right)^{-1}},\end{equation}  we define the
associated mapping $\zeta(Z,W)$ for two arbitrary commuting matrices
$Z$ and $W$ just by replacing $z$ by $Z$ and $w$ by $W$ whenever the
denominator yields a regular matrix. Sometimes we are given the
rational mapping $\zeta(Z,W)$ first and then we define the rational
complex function just by replacing the matrices $Z$ and $W$ by the
complex variables $z$ and $w$, respectively. The above definitions
are straightforward extended to functions and mappings of more than
two complex variables.

We will be mainly concerned with the case in which $z=\tau J$ and $w= \tau P$, where $J$ and $P$
are defined in  (\ref{lin-system}) and  (\ref{matrixP}), respectively. It should be noticed that for instance
the $(i,j)$-element of the matrix $M_\nu$, see (\ref{Mnu-Qnu}), would be given by (observe that it is a matrix itself)
$$ M_{ij}(\tau J, \tau P)= (e_i^T\otimes I) (I_s\otimes I_m- T_\nu\otimes \tau P)^{-1}
(A\otimes \tau J - T_\nu \otimes \tau P)(e_j\otimes I),$$
where $e_j$ denotes the $j$-vector of the canonical basis in $\mathbb{R}^s$ and the corresponding
 complex function is $$M_{ij}(z,w)=e_i^T(I_s-wT_\nu)^{-1}(zA-wT_\nu)e_j.$$

Another important point is that despite  of we are considering cases
with a $d$-splitting for $J$ as indicated in (\ref{lin-system}), the
replacement of every $\tau J_j$ by the complex variable $z_j$ and
the definition of
\begin{equation}\label{z-w}
z:=\sum_{k=1}^d z_k, \qquad w:=
\gamma^{-1}\left(1-\prod_{k=1}^d (1-\gamma
z_k)\right),
\end{equation}
simplifies the study to the case of  two complex variables $z$ and
$w$ or well to  the case of mappings acting on the two matrices
$\tau J$ and $\tau P$.

It is worth to mention that our rational mappings and related
complex functions are all well defined whenever Re $z_k \le 0$ for
$k=1,2,\ldots,d$ and $d$ arbitrary, because the existence of the
matrix inverse $(I-wT_\nu)^{-1}$ is guaranteed if and only if
$\displaystyle{(1-\gamma w)^{-1}=\prod_{k=1}^d (1-\gamma z_k)^{-1}}$
exists. It is easily seen the existence of the late expression by
virtue of  $\gamma>0$ and that  all the eigenvalues of the matrices
$J_j,\;(j=1,\ldots,d)$ have a non-positive real part. Moreover, for
any $\nu=1,\ldots,q$ and any $d\ge 1$, we next prove that
\begin{equation}\label{sev-res0a}\begin{array}{c}
\displaystyle{\sup_{\mbox{\scriptsize Re}\:z_k\:\le 0,
\;k=1,\ldots,d} |Q_\nu(z,w)|< + \infty, \quad
\sup_{\mbox{\scriptsize Re}\:z_k\:\le 0, \;k=1,\ldots,d}
|M_\nu(z,w)|< + \infty,}\\ z \;\mbox{\rm and } w\; \mbox{\rm defined
in (\ref{z-w}).} \end{array}
\end{equation} This see this, observe that $(T_\nu-\gamma I)$ is a
nilpotent matrix fulfilling $(T_\nu-\gamma I)^s=0$ and that
$$ \begin{array}{rcl} Q_\nu(z,w)&=& (I-w T_\nu)^{-1}=\left((1-w\gamma) I-w (T_\nu-\gamma I)\right)^{-1}\\
&=& (1-w\gamma)^{-1}\left( I-\frac{w}{1-w \gamma} (T_\nu-\gamma
I)\right)^{-1}=(1-w\gamma)^{-1}\sum_{j=0}^{s-1}\left(\frac{w}{1-w\gamma}\right)^j
(T_\nu-\gamma I)^{j}
\end{array}
$$
and
$$ \begin{array}{rcl} M_\nu(z,w)&=& Q_\nu(z,w)(z A-w T_\nu)= \frac{z}{1-w\gamma}\left(\sum_{j=0}^{s-1}
(\frac{w}{1-w\gamma})^j (T_\nu-\gamma I)^{j}\right)A \\[0.2pc] &-&
\frac{w}{1-w\gamma}\left(\sum_{j=0}^{s-1}(\frac{w}{1-w\gamma})^j
(T_\nu-\gamma I)^{j}\right)T_\nu.
\end{array}
$$
Hence the boundedness of $Q_\nu(z,w)$ and $M_\nu(z,w)$ follows from
the boundedness of
$$ \begin{array}{rcl}  |\frac{1}{1-w\gamma}|&=&|\prod_{k=1}^d (1-\gamma
z_k)^{-1}| \le 1,  \\
|\frac{w}{1-w\gamma}|&=&\gamma^{-1}|1-\frac{1}{1-w\gamma}| \le
\gamma^{-1}(1+1)=2\gamma^{-1},
\end{array}
$$
and from the next lemma.\hfill $\Box$ }
\end{remark}
\begin{lemma}\label{sev-lema-0} For any $d=2,3,\ldots$, and
$z$ and  $w$ defined in (\ref{z-w}), we have that
$$\begin{array}{c}
\displaystyle{\sup_{\mbox{\scriptsize Re}\:z_k\:\le 0 \atop
k=1,\ldots,d} \left|\frac{z}{1-\gamma w}\right|=\gamma^{-1}\left(
\frac{(d-1)^{d-1}}{d^{d-2}}\right)^{1/2}.}
\end{array}
$$
\end{lemma}

\noindent {\bf Proof.} The third equality below follows from the
Maximum Modulus principle, which says that the Maximum Modulus is
reached at the boundary of the open region for complex analytical
functions,
$$\begin{array}{rcl}
\displaystyle{\sup_{\mbox{\scriptsize Re}\:z_k\:\le 0 \atop
k=1,\ldots,d} \left|\frac{z}{1-\gamma w}\right|}&=&
\gamma^{-1}\displaystyle{\sup_{\mbox{\scriptsize Re}\:z_k\:\le 0
\atop k=1,\ldots,d} \left|\frac{\gamma z}{1-\gamma w}\right|}=
\gamma^{-1}\displaystyle{\sup_{\mbox{\scriptsize Re}\:u_k\:\le 0
\atop k=1,\ldots,d}
\left|\frac{u_1+u_2+\ldots+u_d}{\prod_{k=1}^d(1-u_k)}\right|}\\[0.5pc]
&=&
\gamma^{-1}\displaystyle{\left|\frac{(y_1+y_2+\ldots+y_d)i}{\prod_{k=1}^d\sqrt{1+(y_k)^2}}\right|}=
\gamma^{-1}\displaystyle{\max_{x_k\:\ge 0 \atop
k=1,\ldots,d}\left(\frac{(x_1+x_2+\ldots+x_d)^2}{\prod_{k=1}^d(1+(x_k)^2)}\right)^{1/2}.}
\end{array}
$$
The computation of the extrema  by making zero  the gradient of the
real function of several variables ($x_1,\ldots,x_d$)  gives the
maximum for $x_1=x_2=\ldots=x_d=(d-1)^{-1/2}.$ The proof follows
after substituting above this value. \hfill $\Box$

\begin{definition}\label{def-estability} {\rm A method of the form (\ref{AMF-RK1}) is said to be
$A(\theta)$-stable for a $d$-splitting, if and only if
$$
 |R_q(z,w)| \le 1,\quad \forall z,w  \; \mbox{\rm given by (\ref{z-w}) whenever }
 z_k \in \mathcal{W}(\theta),\;k=1,2,\ldots,d,
$$
where (we consider that the argument of a no-null complex number
ranges in $[-\pi,\pi)$)
\begin{equation}\label{ss1-1}
 \mathcal{W}(\theta):= \{ u \in \mathbb{C}: u=0 \; \mbox{or}\;
|\mbox{arg}(-u)| \le \theta \}.\end{equation}}
\end{definition}

\subsection{Analysis of the Local  Errors}

Next, we study the {\it time-space local  errors} $l_n$ given by
(\ref{local-errors}).   We will see that the time-space local error
$l_n$ is composed of two terms, $l_n^{[2]}$ related to the predictor
used in the {\sf AMF$_q$-RK} method and $l_n^{[1]}$ related to  the
quadrature associated with the underlying Runge-Kutta method.

\begin{lemma}\label{lem-loc-err}
If the linear system has continuous derivatives $u_h^{(k)}(t)$ up to
order $p+1$ in $[0,t^*]$ and the underlying RK method has stage
order $\ell\ge 1$ ($\ell \le p$), i.e.
$$Ac^{j-1}=j^{-1}c^j,\qquad  b^Tc^{j-1}=j^{-1},\quad j=1,2,\ldots,\ell.$$ Then, the local
error $l_n$ in (\ref{local-errors}) of the {\sf AMF$_q$-RK} method
is given by
\begin{equation}\label{local-errors1}
\begin{array}{lll} l_n&=&l_n^{[1]}+
l_n^{[2]},\\ l_n^{[1]}&:=&(\ss^T\otimes I)\left( Q_q + \sum_{j=q}^1
(\prod_{i=q}^{j} M_i)Q_{j-1}-\prod_{i=q}^{1} M_i \right) \hat{D}_n + \delta_n,\\
l_n^{[2]}&:=&(\ss^T\otimes I) (\prod_{i=q}^{1} M_i)\: \triangle
u_h(t_n),
\end{array}
\end{equation}
with
\begin{equation}\label{triangle} \begin{array}{rcl}\triangle
u_h(t_n)&:=&(u_h(t_n+c_i\tau)-u_h(t_n))_{i=1}^s=\sum_{j=1}^p\frac{\tau^j}{j!}(c^j\otimes
I) u^{(j)}_h(t_n)\\[0.3pc] & + & \displaystyle{\frac{\tau^{p+1}}{p!}\left(\int_0^1 (c_i-\theta)_+^p
 u^{(p+1)}_h(t_n+\theta \tau) d\theta\right)_{i=1}^s.}
\end{array}
\end{equation}
\end{lemma}
and (we use, $(x)_+:=x$ if $x\ge 0$ and $(x)_+:=0$ otherwise)
\begin{equation}\label{resid-2} \begin{array}{rcl}
\hat{D}_n &=& \displaystyle{\sum_{j=\ell+1}^p \frac{\tau^j}{j!}
\left((c^j-jAc^{j-1})\otimes u^{(j)}_h(t_n) \right)+ \tau^{p+1}
\int_0^1 \left(\varphi(\theta)
\otimes u^{(p+1)}_h(t_n+ \theta \tau)\right) d\theta + }\\[0.5pc] & & \tau
(A\otimes I) \left(\sigma_h(t_n+c_i\tau)\right)_{i=1}^s;\quad
\varphi(\theta)=\displaystyle{\frac{1}{p!}\left((c_i-\theta)_+^p-p\sum_{j=1}^s
a_{ij}(c_j-\theta)_+^{p-1}\right)_{i=1}^s} \\[0.8pc]
\delta_n &=&  \displaystyle{\sum_{j=\ell+1}^p \frac{\tau^j}{j!}
(1-\ss^T c^j) u^{(j)}_h(t_n)+ \tau^{p+1} \int_0^1\phi(\theta)\:
u^{(p+1)}_h(t_n+ \theta
\tau) d\theta},\\[0.3pc]
\phi(\theta)&=&\displaystyle{\frac{1}{p!}\left((1-\theta)^p-\sum_{j=1}^s
\ss_{j}(c_j-\theta)_+^{p}\right)}.
\end{array}
\end{equation}

\noindent {\bf Proof.} Let us define
\begin{equation}\label{sev-res1}
\hat{D}_n:= (u_h(t_n+c_i\tau))_{i=1}^s-e\otimes u_h(t_n) - \tau
(A\otimes I) (f_h(t_n+c_i
\tau,u_h(t_n+c_i\tau))_{i=1}^s.\end{equation} From
(\ref{spatialerrors}), it follows that
\begin{equation}\label{sev-res1a}\hat{D}_n=
(u_h(t_n+c_i\tau))_{i=1}^s-(u_h(t_n))_{i=1}^s - \tau (A\otimes I)
(u'_h(t_n+c_i\tau)-\sigma_h(t_n+c_i\tau))_{i=1}^s.\end{equation}
Now, by using the Taylor expansion with integral remainder (below
$\zeta(x)$ denotes a generic function having $r+1$-continuous
derivatives in an adequate interval)
\begin{equation}\label{sev-res2}
\displaystyle{\zeta(t_n+x)=\sum_{l=0}^r \frac{x^l}{l!}
\zeta^{(l)}(t_n) + \frac{x^{r+1}}{r!}\int_0^1(1-\theta)^r
\zeta^{(r+1)}(t_n+\theta x) d\theta},
\end{equation}
and applying it conveniently to $u_h(t_n+c_i\tau)$ and
$u'_h(t_n+c_i\tau)$ in (\ref{sev-res1a}) with $r=p$ and $r=p-1$
respectively, we deduce after some computations, the  expression for
$\hat{D}_n$ in (\ref{resid-2}). Observe that order stage $\ell$ for
the Runge-Kutta method implies that $c^j-jAc^{j-1}=0,
\;\ss^Tc^j-1=0,\;j=1,\ldots,\ell$. The expression for $\delta_n$ is
obtained in a similar way, but taking into account that this time we
define,
\begin{equation}\label{sev-res3} \delta_n:= u_h(t_n+\tau)- \varpi u_h(t_n)
- \sum_{j=1}^s \ss_j u_h(t_n+c_j\tau).\end{equation} Let us now take
$\hat{U}_n:=(u_h(t_n+c_i\tau))_{i=1}^s$ and
$\Delta_n^\nu:=\hat{U}_n- U_n^\nu,$ where $ U_n^\nu$ are the
iterates obtained by the scheme (\ref{AMF-RK1}) when the predictor
$U_n^0=e\otimes u_h(t_n)$ is taken on the exact solution of the PDE
at $t_n$, i.e. $y_n=u_h(t_n)$. This gives as  solution, see
(\ref{AMF-RK1})
\begin{equation}\label{sev-res4}
y_{n+1}= \varpi u_h(t_n) + (\ss^T\otimes I)U_n^q.
\end{equation}
From (\ref{sev-res3}) and (\ref{sev-res4}) it follows
\begin{equation}\label{sev-res5}
l_n=u_h(t_{n+1})-y_{n+1}=   (\ss^T\otimes I)\Delta_n^q + \delta_n.
\end{equation}
In order to compute $\Delta_n^q$ we  insert the expression for
$U_n^\nu$ in (\ref{AMF-RK1}).  It follows for the semi-linear
problem (\ref{lin-system}) that
$$
\begin{array}{lll}
(I\otimes I- T_\nu\otimes \tau P)(\Delta_n^\nu-\Delta_n^{\nu-1})&=&
-D(t_n,\tau,u_h(t_n),U_n^{\nu-1})\\&=&-(I\otimes I- A\otimes \tau J)
\Delta_n^{\nu-1}+ \hat{D}_n,\end{array} \quad (\nu=1,2,\dots,q).
$$
This implies that $ \Delta_n^\nu = M_\nu \Delta_n^{\nu-1}+ Q_\nu
\hat{D}_n$, $1\leq \nu\leq q$, and from this recurrence
$$
\Delta_n^q=\displaystyle{\left( Q_q + \sum_{j=q}^1 (\prod_{i=q}^{j}
M_i)Q_{j-1} - \prod_{i=q}^{1} M_i\right) \hat{D}_n+ (\prod_{i=q}^{1}
M_i)\: \Delta_n^0},$$ with $\Delta_n^0=\triangle u_h(t_n)$ in
(\ref{triangle}). Now, from this expression and from
(\ref{sev-res5}) the formula (\ref{local-errors1}) follows. \hfill
$\Box$

\begin{theorem}\label{th3-0}
Consider a family of matrices $\{J_{k,h}\}_{k=1}^d$   and $P_h$,
$h\rightarrow 0^+$, as   given  in (\ref{lin-system}) and
(\ref{matrixP}), respectively. Assume that (\ref{jordan}) holds and
that
\begin{equation}\label{spect-1}  \bigcup_{k=1}^d \mbox{\rm
Spect}(J_{k,h}) \subseteq \mathcal{W}(\theta),\qquad (h\rightarrow
0^+)\end{equation} is fulfilled for some $\theta\in[0,\pi/2]$. Let
$L(z,w)$ be a complex rational function satisfying
$$\sup_{z_k \in \mathcal{W}(\theta), \;k=1,2,\ldots,d} |L(z,w)| \le 1, \quad
\mbox{\rm $z$ and $w$ given by (\ref{z-w}).} $$
 Then,  we have that
$$\begin{array}{c}
 \Vert L(\tau J, \tau P)^n \Vert \le C^*, \quad 0\leq
n\tau\leq t^*, \qquad (\tau,h\rightarrow 0^+).\end{array}$$
\end{theorem}
{\bf Proof.} For simplicity of notations,  we omit the sub-index $h$
in the matrices. By virtue of (\ref{lin-system}), (\ref{jordan}) and
(\ref{sev-res0}) it follows that
%\begin{equation}\label{ss3-0}
$$
\Vert \left(L(\tau J, \tau P)\right)^n \Vert=
\Vert  \Theta \cdot \left(L(\tau \Lambda, \tau
\Upsilon)\right)^n \cdot  \Theta^{-1} \Vert \le C
\Vert \left(L (\tau \Lambda, \tau \Upsilon
)\right)^n\Vert, \quad n\geq 1,
$$
%\end{equation}
where
%\begin{equation}\label{ss3-0a}
$$
\begin{array}{l}
\tau \Lambda:=\sum_{k=1}^d \tau \Lambda_k, \quad
\Lambda_k=\mbox{\rm
Block-Diag}(\Lambda^{(1)}_k,\Lambda^{(2)}_k,\ldots,
\Lambda^{(\vartheta)}_k),
\\ \tau \Upsilon := \gamma^{-1}\left(I-\prod_{k=1}^d (I-\gamma \tau
\Lambda_k)\right).
\end{array}
$$
%\end{equation}
By defining $\tau \Lambda^{(l)}:= \sum_{k=1}^d \tau \Lambda^{(l)}_k$
and $\tau \Upsilon^{(l)}:= \gamma^{-1}\left(I-\prod_{k=1}^d
(I-\gamma \tau \Lambda^{(l)}_k)\right)$, for the  norm considered it
follows that
%\begin{equation}\label{ss3-0b}
$$
\Vert  \left(L(\tau \Lambda, \tau
\Upsilon)\right)^n \Vert = \max_{
l=1,\ldots,\vartheta} \Vert \left(L(\tau
\Lambda^{(l)} , \tau \Upsilon^{(l)})\right)^n
\Vert, \quad n\geq 1.
$$
%\end{equation}
Consider any diagonal block  $\Lambda^{(l)}_k=\lambda_k^{(l)}I + E$
($E\equiv E^{(l)}$ for simplicity of notation. Observe that all the
matrices $E$ are strictly lower triangular and they have  uniform
bounded entries and uniform bounded dimensions, hence  all of them
are nilpotent with nilpotency index $\le N$) and define
$$z_k= \tau \lambda^{(l)}_k,\;1\leq k\leq d,\; z=\sum_{k=1}^d z_k,
\; w= \gamma^{-1}\left(1-\prod_{k=1}^d (1-\gamma z_k)\right),$$ it
follows that,
$$ L(\tau \Lambda^{(l)} , \tau \Upsilon^{(l)})  = L\left(\sum_{k=1}^d
(z_k I+\tau  E),\gamma^{-1}(I-\prod_{k=1}^d (I-\gamma (z_k I+\tau
E))\right).
$$
By defining the function of $d$ complex
variables,
%\begin{equation}\label{ss3-0c}
$$
\psi(w_1,\ldots,w_d):= L\left(\sum_{k=1}^d w_k
,\gamma^{-1}(1-\prod_{k=1}^d (1-\gamma
w_k))\right),
%\end{equation}
$$
we get that $ L(\tau \Lambda^{(l)} , \tau \Upsilon^{(l)})=\psi(z_1
I+\tau E,\ldots,z_d I+\tau E).$ Then, by using the Taylor expansion
for $\psi$ around $\tau=0$ and taking into the nilpotency of the
matrix $E$,
 we
deduce that,
%\begin{equation}\label{ss3-0e}
$$ \begin{array}{l}
\psi(z_1 I+\tau E,\ldots,z_d I+\tau E)= \psi(z_1,\ldots,z_d) I +
\\ \qquad \qquad \displaystyle{ \sum_{l=1}^{N-1}
\frac{\tau^l}{l!} E^l \sum_{i_1+i_2+\ldots+i_d=l}
\frac{\partial^{l}\psi }{\partial^{i_1} z_1\ldots\partial^{i_d} z_d}
(z_1,z_2,\ldots,z_d)}.\end{array}
$$
%\end{equation}
Now, since $L(z,w)\equiv L(z_1,\ldots,z_d)$ and all its partial
derivatives up to order $N$ are uniformly bounded on the wedge
$\mathcal{W}(\theta)$, we can write that
$$ \psi(z_1 I+\tau E,\ldots,z_d I+\tau E) =  \psi(z_1,\ldots,z_d ) I + \tau L^*_{\tau,h},\quad
\Vert L^*_{\tau,h}\Vert \le C^*, \;(\tau,\:h\rightarrow 0^+).$$
From here we get for $0\le
\tau n \le t^*$ that
$$ \Vert \left(\psi(z_1 I+\tau E,\ldots,z_d I+\tau
E)\right)^n \Vert =  \Vert \left (L(z,w) I + \tau
L^*_{\tau,h}\right)^n\Vert  \le (1+\tau C^*)^n \le \exp(t^*C^*).$$
 $\Box$

\subsection{Some mappings and definitions}
For a given mapping $\zeta(X,Y) \in \mathbb{C}^{m,m}$ where $X$ and
$Y$ are two  arbitrary square complex matrices of order $m$ we
define some associated mappings in the following way,
\begin{equation}\label{sev-res7}\begin{array}{c}
\zeta^{[1]}(X,Y):=\left(\zeta(X,Y)-\zeta(X,X)\right)(Y-X)^{-1},
\quad \mbox{\rm
whenever } \det(Y-X)\ne 0,\\
\zeta^{[1]}(X,X):=\lim_{\epsilon\rightarrow 0}
\zeta^{[1]}(X,X+\epsilon I), \quad \mbox{\rm whenever the limit
exists. }
\end{array}
\end{equation}
In a recursive form, when $\det(Y-X)\ne 0$ and $\zeta^{[l]}(X,X)$
exists,  we continue by defining
\begin{equation}\label{sev-res8}\begin{array}{c}
\zeta^{[l+1]}(X,Y):=\left(\zeta^{[l]}(X,Y)-\zeta^{[l]}(X,X)\right)(Y-X)^{-1},
\\
\zeta^{[l+1]}(X,X):=\lim_{\epsilon\rightarrow 0}
\zeta^{[l+1]}(X,X+\epsilon I), \quad l=1,2,\ldots,l^*.
\end{array}
\end{equation}
By assuming  $\det(Y-X)\ne 0$ and the existence of
$\zeta^{[l]}(X,X),\;l=1,2,\ldots,l^*$, it is straightforward to show
by induction that
\begin{equation}\label{sev-res9}
\zeta(X,Y)=\sum_{l=0}^{l^*}\zeta^{[l]}(X,X)(Y-X)^{l} +
\zeta^{[l^*+1]}(X,Y)(Y-X)^{l^*+1}.
\end{equation}
We have considered for  convenience that
$\zeta^{[0]}(X,Y):=\zeta(X,Y).$ It should be noted that the
commutativity of the matrices $X$ and $Y$ is neither necessary in
the definitions above nor in the formula (\ref{sev-res9}).

To have a practical meaning of the mapping $\zeta^{[l]}(X,Y)$ we
show next that assuming  $\zeta(x,y)$ has $l^*$  continuous partial
derivatives regarding the second variable, then it holds that
\begin{equation}\label{sev-res10}
\zeta^{[l]}(X,X)=\frac{1}{l!} \frac{\partial^l \zeta(x,y)}{\partial
y^l} (X,X),\quad l=1,2,\ldots,l^*.
\end{equation}
To see (\ref{sev-res10}), we use the induction.  For $l=0$ it is
true for convenience.  For $l=1$ it is true since
$$ \zeta^{[l]}(X,X)=\lim_{\epsilon \rightarrow 0} \zeta^{[l]}(X,X+\epsilon
I)=\lim_{\epsilon \rightarrow 0} \epsilon^{-1}
\left(\zeta(X,X+\epsilon I)-\zeta(X,X)\right)= \frac{\partial
\zeta}{\partial y} (X,X).$$ Assume it is true up to  $l$, we show it
 for $l+1$ by using (\ref{sev-res9}) in the second equality and the
induction in the third equality below. The L'Hospital formula for
limits (for the indetermination $0/0$) is used $l+1$ times in the
fourth equality,
$$\begin{array}{rcl}
\zeta^{[l+1]}(X,X)&=&\displaystyle{\lim_{\epsilon\rightarrow 0}
\zeta^{[l+1]}(X,X+\epsilon I)=\lim_{\epsilon\rightarrow
0}\frac{\zeta(X,X+\epsilon
I)-\sum_{j=0}^{l}\zeta^{[j]}(X,X)(\epsilon I)^{j}}{(\epsilon
I)^{l+1}}}\\[0.5pc] &= & \displaystyle{\lim_{\epsilon\rightarrow 0}
\frac{\zeta(X,X+\epsilon I)-\sum_{j=0}^{l}\frac{\epsilon^j}{j!}
\frac{\partial^j \zeta(x,y)}{\partial y^l}(X,X)}{\epsilon^{l+1}}}\\[0.5pc]
& = &\displaystyle{\lim_{\epsilon\rightarrow 0}  \frac{1}{(l+1)!}
\frac{\partial^{l+1} \zeta(x,y)}{\partial y^{l+1}}(X,X+\epsilon I) =
\frac{1}{(l+1)!} \frac{\partial^{l+1} \zeta}{\partial y^{l+1}}(X,X)}
\end{array}$$

These results can be trivially extended to  vectors (and matrices),
namely $(\zeta_{ij}(X,Y))\in \mathbb{C}^{q_1m,q_2m}$, by applying
them to each component $\zeta_{ij}(X,Y)\in \mathbb{C}^{m,m}$.
Sometimes we will make use of this kind of vectors as we will see in
the next section.

\subsection{Bounds for the local errors}

 The forthcoming convergence results for {\sf AMF$_q$-RK} methods
 are based in the Lemma II.2.3 \cite[p. 162]{HV},
which can be stated as follows
\begin{lemma}\label{sev-lema-glob-err}
Assume that the global errors $\epsilon_n\equiv\epsilon_n(\tau;h)$,
of a one-step method satisfy the recursion (\ref{global-errors1}),
where the local errors $l_n$ can be split (uniformly on $h$ and
$\tau$) as
\begin{equation}\label{ln-domi}
l_n=\left(R_q(\tau J,\tau P)-I\right) \phi(t_n)\tau^\mu h^\alpha +
\tau \mathcal{O}(\tau^\nu h^\beta),\quad n=0,1,\ldots, t^*/\tau-1,
\end{equation}
where the function $\phi(t)$ and its first derivative regarding $t$
are uniformly bounded, then the stability condition
\begin{equation}\label{sev-estab}
\sup_{1\le n\le t^*/\tau \atop \tau\rightarrow 0^+, \;h\rightarrow
0^+} \Vert R_q(\tau J,\tau P)^n \Vert\le C,
\end{equation}
 implies that the global errors uniformly fulfil
\begin{equation}\label{sev-global-err} \epsilon_{n} = \mathcal{O}(\tau^\mu h^\alpha) +
\mathcal{O}(\tau^\nu h^\beta),\quad n=1,\ldots, t^*/\tau, \quad
\tau\rightarrow 0^+, h\rightarrow 0^+.\end{equation}
\end{lemma}

{\sf General Assumptions on the semilinear problem.}

{\it To bound the local errors and consequently  the global errors
we henceforth assume that the exact PDE solution $u_h(t)$ confined
to the spatial grid and the semilinear problem (\ref{lin-system})
 fulfil (\ref{spatialerrors})-(\ref{norm}),
(\ref{jordan}) and (\ref{spect-1}) for some $\theta\in [0,\pi/2]$,
and that the following hypotheses (related the matrices $J$ and $P$)
hold for some constants (not necessarily positive) $\alpha_l,\;
\beta_l $ and $\eta$ and some nonnegative integer $l^*$, whenever
$h\rightarrow 0^+$ and $\tau\rightarrow 0^+$,
\begin{equation}\label{P-hypo}
\begin{array}{rcl} {\bf (P1)} & & \left\{ \begin{array}{l}
 (P-J)^lu_h^{(k)}(t)=\tau^l h^{\alpha_l} \:
\mathcal{O}(1), \\
  (P-J)^{l+1} u_h^{(k)}(t)=\tau^{l+1} h^{\beta_{l+1}} J \:
\mathcal{O}(1) \end{array} \right\} \quad  {l=0,1,\ldots,l^* \atop
 k=1,2,\ldots,p+1.} \\[0.5pc]
{\bf (P2)} & &  J^\eta u_h^{(k)}(t)=\mathcal{O}(1), \quad
k=1,2,\ldots,p+1, \;\mbox{\rm for some } \eta.
\end{array}
\end{equation}
It should be noticed that always $\alpha_0=0$, because the
derivatives (up to some order) of the exact solution are uniformly
bounded, i.e. $u_h^{(k)}(t)=\mathcal{O}(1), \;t\in
[0,t^*],\;k=0,1,\ldots,p+1$.}

\begin{theorem}\label{sev-th-1} Assume  that the
Runge-Kutta method has stage order $\ell$ and that
\begin{equation}\label{sev-equ-1}\sup_{z_k\in \mathcal{W}(\theta),\atop k=1,2,\ldots,d} |z/(R_q(z,w)-1)|\le C, \quad z\;\mbox{\rm and } w \; \mbox{\rm
given by } (\ref{z-w}). \end{equation} Then for the {\sf AMF$_q$-RK}
method we have that,
$$l_n^{[1]}= \mathcal{O}(\tau h^r) + \mathcal{O}(\tau^{\ell+1}),
\qquad (\tau\rightarrow 0, \; h\rightarrow 0),$$ and
$$l_n^{[1]}= \tau h^r\mathcal{O}(1) + \tau^{\ell+1}(R_q(\tau J, \tau P)-I)\left(\mathcal{O}(1) +
\tau h^{\beta_1} \mathcal{O}(1)\right), \; \tau\rightarrow 0, \;
h\rightarrow 0.
$$
\end{theorem}

{\bf Proof.} According to Lemma \ref{lem-loc-err} the term
$l_n^{[1]}$ of the local error is given by,
\begin{equation}\label{sev-res11}
l_n^{[1]}= \xi(\tau J, \tau P) \hat{D}_n + \delta_n,
\end{equation}
where
\begin{equation}\label{sev-res12}
\xi(\tau J, \tau P):=(\ss^T\otimes I)\left( Q_q(\tau J, \tau P) +
\sum_{j=q}^1 (\prod_{i=q}^{j} M_i(\tau J, \tau P))Q_{j-1}(\tau J,
\tau P)-\prod_{i=q}^{1} M_i(\tau J, \tau P) \right)
\end{equation}
From Remark \ref{sev-remark-0} we have that ($e_j$ denotes the
$j$-vector of the canonical basis)
$$ \sup_{\mbox {\tiny Re } z_k \le 0\atop k=1,\ldots,d}|\xi(z,w)e_j|\le C, \quad (j=1,\ldots,s),\quad
z,\; w \;\mbox{\rm given by (\ref{z-w})}. $$ From Theorem
\ref{th3-0} this implies  that
$$ \max_{j=1,\ldots,s} \Vert \xi(\tau J, \tau P)(e_j\otimes I) \Vert \le C', \quad
\tau \rightarrow 0^+, \quad h \rightarrow 0^+.$$ Then, from
(\ref{resid-2}) in Lemma \ref{lem-loc-err} the first bound for
$l_n^{[1]}$ follows.

For the second bound, we separate in (\ref{sev-res11}) the
$\tau^{\ell+1}$-term from the others, take into account
(\ref{sev-res12}) and Lemma \ref{lem-loc-err}, we get
\begin{equation}\label{sev-res13}\begin{array}{rcl}
l_n^{[1]}&=& \frac{\tau^{\ell+1}}{(\ell+1)!} \left(\xi(\tau J, \tau
P)\left((c^{\ell+1}-(\ell+1)Ac^\ell)\otimes I \right) + (1-\ss^T
c^{\ell+1}) I
 \right) u_h^{(\ell+1)}(t_n) + \circledR,\\[0.3pc]
\mbox{\rm where }& &\circledR= \mathcal{O}(\tau^{\ell+2}) +
\mathcal{O}(\tau h^r).
\end{array}
\end{equation}
Next, we define the mapping (assume that $J$ is regular only to
simplify the proof)
\begin{equation}\label{sev-res14}
\upsilon(\tau J, \tau P):=(R_q(\tau J, \tau P) - I)^{-1}\left(
\xi(\tau J, \tau P)\left((c^{\ell+1}-(\ell+1)Ac^\ell)\otimes
I\right) + (1-\ss^T c^{\ell+1}) I\right).
\end{equation}
By using the assumption (\ref{sev-equ-1}), the bounds in Remark
\ref{sev-remark-0} and Lemma \ref{sev-lema-0}, it is not very
difficult to  see that
\begin{equation}\label{sev-res14a} \begin{array}{c}
\displaystyle{\sup_ {z \in \mathcal{W}(\theta)}|\upsilon(z, z)|<
+\infty. \quad \sup_{z_k \in \mathcal{W}(\theta) \atop
k=1,2,\ldots,d }|z\upsilon^{[1]}(z, w)|< +\infty, \; \mbox{\rm $z$
and $w$ given by (\ref{z-w}).}}\end{array}
\end{equation} Then, from (\ref{sev-res13}) it follows that,
\begin{equation}\label{sev-res15}\begin{array}{rcl}
l_n^{[1]}&=& \frac{\tau^{\ell+1}}{(\ell+1)!} \left(R(\tau J, \tau
P)-I\right)\upsilon(\tau J, \tau P) u_h^{(\ell+1)}(t_n)  +
\circledR \\[0.3pc] &=& \frac{\tau^{\ell+1}}{(\ell+1)!} \left(R(\tau J, \tau
P)-I\right)\left(\upsilon(\tau J, \tau J)+ \upsilon^{[1]}(\tau J,
\tau P)(\tau P-\tau J)\right)u_h^{(\ell+1)}(t_n)  +
\circledR  \\[0.3pc] &=& \circledR + \frac{\tau^{\ell+1}}{(\ell+1)!} \left(R(\tau J, \tau
P)-I\right)\upsilon(\tau J, \tau J)u_h^{(\ell+1)}(t_n)  \\[0.3pc] &+& \frac{\tau^{\ell+1}}
{(\ell+1)!} \left(R(\tau J, \tau
P)-I\right)\upsilon^{[1]}(\tau J, \tau P)(\tau J)(J^{-1}(P-J))
u_h^{(\ell+1)}(t_n)\\[0.3pc] &=& \circledR + \frac{\tau^{\ell+1}}{(\ell+1)!}
\left(R(\tau J, \tau P)-I\right) \mathcal{O}(1) +
\frac{\tau^{\ell+1}}{(\ell+1)!}\left(R(\tau J, \tau P)-I\right)
\mathcal{O}(\tau h^{\beta_1})\quad \mbox{\rm \hfill $\Box$}
\end{array}
\end{equation}

For the analysis of the local error term $l_n^{[2]}$ in
(\ref{local-errors1}), we define the mappings
\begin{equation} \begin{array}{lll}\label{H-eq}
\psi_{q}(\tau J,\tau X)&:=& (\ss^T\otimes I)\prod_{j=q}^1 M_j(\tau
J,\tau X) \in \mathbb{C}^{m,sm}, \\[0.3pc] \zeta_{q}(\tau J,\tau X)&:=& \left(R_q(\tau J,\tau
X)-I\right)^{-1} \psi_{q}(\tau J,\tau X) \in \mathbb{C}^{m,sm},
\end{array}
\end{equation}
and their associated  vector complex  functions
\begin{equation}
\begin{array}{lll}\label{H-eq1} \psi_{q}(z,w)&:=&
\ss^T\prod_{j=q}^1 (I-wT_j)^{-1}(zA-wT_j) \in \mathbb{C}^{1,s}, \\[0.3pc]
\zeta_{q}(z,w)&:=& \left(R_q(z,w)-1\right)^{-1} \psi_{q}(z,w) \in
\mathbb{C}^{1,s}.
\end{array}
\end{equation}
These mappings will play a mayor role in the proof of the
convergence results. It must be remarked that whereas $\Vert
\psi_{q}(z,w)\Vert_2$ is uniformly bounded when $z$ and $w$ are
given by (\ref{z-w}), the vector
$\zeta_{q}(z,w)=\mathcal{O}(z^{-1})$ as $z\rightarrow 0$ due to the
fact that (see (\ref{estab-1}))
\begin{equation}\label{sev-eqq2}
R_q(z,w)-1= \ss^T\left(Q_q(z,w) + \sum_{j=q}^1 (\prod_{i=q}^{j}
M_i(z,w))Q_{j-1}(z,w)-\prod_{i=q}^{1} M_i(z,w)\right) c z.
\end{equation}
Hence $\zeta_{q}(z,w)$  is not bounded in general for $z$ and $w$
 given by (\ref{z-w}). However, $\zeta_{q}(z,z)$ is uniformly
 bounded as long as $R_q(z,z)-1\ne 0$ for $z\in \mathcal{W}(\theta)
 \backslash \{0\}$.

 From (\ref{local-errors1}), by using (\ref{sev-res9}), we
deduce that,
\begin{equation}\label{sev-res16} \begin{array}{rcl}l_n^{[2]}&=& (R_q(\tau J, \tau P)-I)
\sum_{j=0}^{l^*} \zeta_q^{[j]}(\tau J,
\tau J) (I_s\otimes (\tau(P-J))^j)\Delta_h(t_n)\\[0.3pc] & + & (R_q(\tau J, \tau
P)-I) \zeta_q^{[l^*+1]}(\tau J, \tau P)(I_s\otimes
(\tau(P-J))^{l^*+1})\Delta_h(t_n).\end{array}
\end{equation}

 Next, we provide some convergence results for different
kind of {\sf AMF$q$-RK} methods, which depends on the Runge-Kutta
method on which the {\sf AMF$_q$-RK} is based on. We start  with
Theorem \ref{sev-th-2} that meets applications for  DIRK methods
(Diagonally Implicit Runge-Kutta) and SIRK methods (Single Implicit
Runge-Kutta) and then with  Theorems \ref{sev-th-3}, \ref{sev-th-4}
and \ref{sev-th-5} which meet applications in the {\sf AMF$q$-Rad}
methods presented in section two. Of course, the assumptions {\bf
(P1)-(P2)} will be always assumed for some integers $l^*\ge 0$,
$\ell\ge 1, \:p\ge 1$.

\begin{theorem}\label{sev-th-2} If $T_\nu=A, \;\nu=1,\ldots,q$, with the Runge-Kutta coefficient matrix $A$
having unique eigenvalue $\gamma>0$ (with multiplicity $s$), then
the local errors ($l_n=l_n^{[1]}+ l_n^{[2]}$) fulfil
$$\left. \begin{array}{rcl}
l_n^{[1]}&=& \mathcal{O}(\tau h^r) + \tau^{\ell+1}(R(\tau J,\tau
P)-I) (\mathcal{O}(1)+ \mathcal{O}(\tau h^{\beta_1})),
\\
l_n^{[2]}&=& \tau^{2l+2}h^{\beta_{l+1}}(R(\tau J,\tau P)-I)
\mathcal{O}(1),\; l=0,1,\ldots,\tilde{l}, \\
\tilde{l}&=&\max\{0,\min\{q-2,l^*\}\}.
\end{array} \right\} \; (\tau\rightarrow 0^+, \; h\rightarrow 0^+). $$
If the method is A$(\theta)$-stable for a $d$-splitting and
(\ref{spect-1}) holds, then for any $l=0,1,\ldots,\tilde{l}$,
 the global errors fulfil (whenever $\tau\rightarrow 0^+$ and $h\rightarrow 0^+$) that,
$$\epsilon_{n,h} = \mathcal{O}( h^r)+
\tau^{\ell}\min\{1,\max\{\tau,\tau^2 h^{\beta_1}\}\} \mathcal{O}(1)
+\mathcal{O}(\tau^{2l+2} h^{\beta_{l+1}}); \;
n=1,2,\ldots,t^*/\tau.$$
\end{theorem}

{\bf Proof.} The expression of $l_n^{[1]}$ was seen in Theorem
\ref{sev-th-1}. In order to show the expression for $l_n^{[2]}$, we
start  by deducing  from (\ref{H-eq}) and  (\ref{estab-1}) that
\begin{equation}\label{sev-eqq-1} \begin{array}{rcl}
\zeta_q(z,w)&=&(R_q(z,w)-1)^{-1} \ss^T\left((I-wA)^{-1}A\right)^q
(z-w)^q,\\[0.3pc]
R_q(z,w)-1&=&\ss^T\left(\left((I-wA)^{-1}A\right)^q
(z-w)^q-I\right)(zA-I)c z.
\end{array}
\end{equation}
From (\ref{sev-res10}) we have that
$\displaystyle{\zeta_q^{[l]}(z,z)=\frac{1}{l!}\frac{\partial^l
\zeta_q}{\partial w^l}}(z,z)$. From here and from (\ref{sev-eqq-1})
it follows that
$$\zeta_q^{[l]}(z,z)=0,\quad l=0,1,\ldots,\tilde{l}.$$
From (\ref{sev-res16}) by taking $\tilde{l}$ as upper index, for any
$l=0,1,\ldots,\tilde{l}$, we have that
$$ \begin{array}{rcl}
\begin{array}{rcl}l_n^{[2]}&=& (R_q(\tau J, \tau
P)-I) \zeta_q^{[l+1]}(\tau J, \tau P)(I_s\otimes
(\tau(P-J))^{l+1}\Delta_h(t_n) \\[0.3pc] &=& \tau^{l}(R_q(\tau J, \tau
P)-I) \left( \zeta_q^{[l+1]}(\tau J, \tau P)(I_s\otimes \tau
J)\right)(I_s\otimes J^{-1}(P-J)^{l+1})
 \Delta_h(t_n) \\[0.3pc] &=& \tau^{l}(R(\tau J,\tau P)-I)
\mathcal{O}(1) (I_s\otimes J^{-1}(P-J)^{l+1})(\tau c\otimes
u_h'(t_n) + \tau^2\mathcal{O}(1))\\[0.3pc]
&=&\tau^{2l+2}h^{\beta_{l+1}}(R(\tau J,\tau P)-I) \mathcal{O}(1).
\end{array}
\end{array} $$
To see the bound for the global errors we apply  Lemma
\ref{sev-lema-glob-err}. The bounds for the local errors $l_n$ have
been obtained above (see also Theorem \ref{sev-th-1} for
$l_n^{[1]}$). The boundedness of the powers of $R_q(\tau J, \tau P)$
as indicated in (\ref{sev-estab}) follows from Theorem \ref{th3-0}
by taking into account the A$(\theta)$-stability of the method for
the $d$-splitting and that (\ref{spect-1}) holds. Now from Lemma
\ref{sev-lema-glob-err} the proof is accomplished. \hfill $\Box$

\begin{theorem}\label{sev-th-3} For AMF$_q$-RK methods with $\gamma>0$ and
satisfying $(A-T_1)c=0$, we have that
$$
l_n^{[2]}= \tau^{2}(R(\tau J,\tau P)-I)\left( \mathcal{O}(1) +
h^{\beta_1}\mathcal{O}(1)\right), \; (\tau\rightarrow 0^+, \;
h\rightarrow 0^+). $$ Additionally  if the method is
A$(\theta)$-stable for a $d$-splitting and (\ref{spect-1}) holds,
then for $\tau\rightarrow 0^+$ and $h\rightarrow 0^+$, the global
errors fulfil
$$\epsilon_{n,h} = \mathcal{O}(
h^r)+  \tau^{\ell}\min\{1,\max\{\tau,\tau^2 h^{\beta_1}\}\}
\mathcal{O}(1) +\tau^{2}\left(\mathcal{O}(1)+
h^{\beta_1}\mathcal{O}(1) \right); \; n=1,2,\ldots,t^*/\tau.$$
\end{theorem}

{\bf Proof.} The expression of $l_n^{[1]}$ was seen in Theorem
\ref{sev-th-1}. In order to show the expression for $l_n^{[2]}$,
from (\ref{sev-res16}) by setting $l^*=0$ we get that (observe that
$\zeta_q(z,z)c=0$ because $(A-T_1)c=0$. This expression is used in
the third equality below)
$$ \begin{array}{rcl}
\begin{array}{rcl}l_n^{[2]}&=& (R_q(\tau J, \tau
P)-I) \zeta_q(\tau J, \tau J)\Delta_h(t_n) \\[0.3pc] &+&(R_q(\tau J, \tau
P)-I) \zeta_q^{[1]}(\tau J, \tau P)(I_s\otimes
(\tau(P-J)))\Delta_h(t_n) \\[0.3pc] &=&
(R_q(\tau J, \tau
P)-I) \zeta_q(\tau J, \tau J)\left(\tau c \otimes I + \tau^2 \mathcal{O}(1)\right) \\[0.3pc] &+&
(R_q(\tau J, \tau P)-I) \left(\zeta_q^{[1]}(\tau J, \tau P)(I\otimes
\tau J)\right)(I_s\otimes
(J^{-1}(P-J)))\left(\tau  \mathcal{O}(1)\right) \\[0.3pc] &=&
(R(\tau J,\tau P)-I) \left(\tau^{2} \mathcal{O}(1)\right) + (R(\tau
J,\tau P)-I)
\mathcal{O}(1)\:\left(\tau^{2}h^{\beta_1}\mathcal{O}(1)\right).
\end{array}
\end{array} $$
This provides the bound for the local errors $l_n^{[2]}$. The
boundedness of the powers of $R_q(\tau J, \tau P)$ as indicated in
(\ref{sev-estab}) follows from Theorem \ref{th3-0} by taking into
account the A$(\theta)$-stability of the method for the
$d$-splitting and that (\ref{spect-1}) holds. Now, from the bounds
for the local error and from  Lemma \ref{sev-lema-glob-err} the
proof follows. \hfill $\Box$

\begin{theorem}\label{sev-th-4} For AMF$_q$-RK methods with $\gamma>0$ and
satisfying $$\sup_{\mbox{\tiny Re$\:z$} \:\le\: 0, \:z\ne 0}\Vert
z^{-\eta} \zeta_q(z,z)\Vert_2 < +\infty,$$ with $\eta$ given  in
{\bf (P2)} we have that
$$
l_n^{[2]}= (R(\tau J,\tau P)-I)\left( \mathcal{O}( \tau^{1+\eta}) +
\mathcal{O}( \tau^2 h^{\beta_1})\right), \; (\tau\rightarrow 0^+, \;
h\rightarrow 0^+). $$ Additionally  if the method is
A$(\theta)$-stable for a $d$-splitting and (\ref{spect-1}) holds,
then for $\tau\rightarrow 0^+$ and $h\rightarrow 0^+$, the global
errors fulfil
$$\epsilon_{n,h} = \mathcal{O}(
h^r)+ \min\{1,\max\{\tau,\tau^2 h^{\beta_1}\}\} \mathcal{O}(
\tau^{\ell}) + \mathcal{O}(\tau^{1+\eta})+ \mathcal{O}(\tau^2
h^{\beta_1}); \; n=1,2,\ldots,t^*/\tau.$$
\end{theorem}

{\bf Proof.} The expression of $l_n^{[1]}$ was seen in Theorem
\ref{sev-th-1}. In order to show the expression for $l_n^{[2]}$,
from (\ref{sev-res16}) by setting $l^*=0$ we get that
$$ \begin{array}{rcl}
\begin{array}{rcl}l_n^{[2]}&=& (R_q(\tau J, \tau
P)-I) \zeta_q(\tau J, \tau J)\Delta_h(t_n) \\[0.3pc] &+&(R_q(\tau J, \tau
P)-I) \zeta_q^{[1]}(\tau J, \tau P)(I_s\otimes
(\tau(P-J)))\Delta_h(t_n) \\[0.3pc] &=&
(R_q(\tau J, \tau
P)-I) \zeta_q(\tau J, \tau J)\left(\tau  \mathcal{O}(1)\right) \\[0.3pc] &+&
(R_q(\tau J, \tau P)-I) \left(\zeta_q^{[1]}(\tau J, \tau P)(I\otimes
\tau J)\right)(I_s\otimes
(J^{-1}(P-J)))\left(\tau  \mathcal{O}(1)\right) \\[0.3pc] &=&
(R_q(\tau J, \tau P)-I) \left(\zeta_q(\tau J, \tau J)(I\otimes (\tau
J)^{-\eta})\right)
\left(I\otimes (\tau J)^{\eta}\right)\left(\tau \mathcal{O}(1)\right) \\[0.3pc] &+&
(R(\tau J,\tau P)-I)
\mathcal{O}(1)\:\left(\tau^{2}h^{\beta_1}\mathcal{O}(1)\right)\\[0.3pc] &=&
R_q(\tau J, \tau P)-I) \left(\mathcal{O}(1) \right)
\left(\tau^{\eta+1} I\otimes J^{\eta} \mathcal{O}(1)\right) \\[0.3pc] &+&
(R(\tau J,\tau P)-I)
\left(\tau^{2}h^{\beta_1}\mathcal{O}(1)\right)\\[0.3pc] &=&
(R(\tau J,\tau P)-I) \left( \mathcal{O}(\tau^{1+\eta}) +
\mathcal{O}(\tau^{2}h^{\beta_1})\right).
\end{array}
\end{array} $$
This provides the bound for the local errors $l_n^{[2]}$. The rest
of the proof follows as in the previous theorems.  \hfill $\Box$

\begin{theorem}\label{sev-th-5} For AMF$_q$-RK methods with $\gamma>0$ and
 $$\begin{array}{c} (A-T_1)c=0,\quad \sup_{\mbox{\tiny Re$\:z$} \:\le\: 0, \:z\ne 0}\Vert z^{-\eta}
\zeta_q(z,z)\Vert_2 < +\infty,\end{array}$$ with $\eta$ given  in
{\bf (P2)} and assuming {\bf (P1)} for $l^*=1$, we have that
$$
l_n^{[2]}= (R(\tau J,\tau P)-I)\left( \mathcal{O}( \tau^{2+\eta}) +
\mathcal{O}( \tau^3 h^{\alpha_1}) +  \mathcal{O}( \tau^4
h^{\beta_2})\right), \; (\tau\rightarrow 0^+, \; h\rightarrow 0^+).
$$ Additionally  if the method is A$(\theta)$-stable for a
$d$-splitting and (\ref{spect-1}) holds, then  the global errors
fulfil
$$\begin{array}{c} \epsilon_{n,h} = \mathcal{O}(
h^r)+ \min\{1,\max\{\tau,\tau^2 h^{\beta_1}\}\} \mathcal{O}(
\tau^{\ell}) + \mathcal{O}(\tau^{2+\eta})+ \mathcal{O}(\tau^3
h^{\alpha_1})+ \mathcal{O}(\tau^4 h^{\beta_2}),\\
n=1,2,\ldots,t^*/\tau,\qquad (\tau\rightarrow 0^+,\; h\rightarrow
0^+).\end{array}$$
\end{theorem}

{\bf Proof.} In order to show the expression for $l_n^{[2]}$, from
(\ref{sev-res16}) by setting $l^*=1$ we get that
$$ \begin{array}{rcl}
\begin{array}{rcl}l_n^{[2]}&=& (R_q(\tau J, \tau
P)-I)\left( \zeta_q(\tau J, \tau J)+ \zeta^{[1]}_q(\tau J, \tau J)(I\otimes \tau(P-J))\right)
\Delta_h(t_n) \\[0.3pc] &+&(R_q(\tau J, \tau
P)-I) \zeta_q^{[2]}(\tau J, \tau P)(I_s\otimes
\tau^2(P-J)^2)\Delta_h(t_n) \\[0.3pc] &=&
(R_q(\tau J, \tau P)-I)\left( \zeta_q(\tau J, \tau J)+
\zeta^{[1]}_q(\tau J, \tau J)(I\otimes \tau(P-J))\right)
\left((\tau c\otimes I)u_h'(t_n)+   \tau^2\mathcal{O}(1)\right)\\[0.3pc] &+&(R_q(\tau J, \tau
P)-I) \zeta_q^{[2]}(\tau J, \tau P)(I_s\otimes
\tau^2(P-J)^2)(\tau \mathcal{O}(1)) \\[0.3pc] &=&
(R_q(\tau J, \tau P)-I)\left( \tau^2 \zeta_q(\tau J, \tau
J)\mathcal{O}(1)+ \tau \zeta^{[1]}_q(\tau J, \tau J)(I\otimes
\tau(P-J)\mathcal{O}(1))\right)
\\[0.3pc] &+&(R_q(\tau J, \tau
P)-I) \left(\zeta_q^{[2]}(\tau J, \tau P)(I_s\otimes \tau
J)\right)(I_s\otimes \tau J^{-1}(P-J)^2)(\tau \mathcal{O}(1))
\\[0.3pc] &=&
(R_q(\tau J, \tau P)-I)\left(\tau^2\left(\zeta_q(\tau J, \tau
J)(\tau J)^{-\eta}\right)(\tau^\eta J^\eta \mathcal{O}(1))+
\mathcal{O}(\tau^3 h^{\alpha_1})\right)
\\[0.3pc] &+&
(R_q(\tau J, \tau P)-I) \left(\mathcal{O}(1)\: \tau^2 J^{-1}(P-J)^2
\mathcal{O}(1)\right)
\\[0.3pc] &=&
(R(\tau J,\tau P)-I)\left( \mathcal{O}( \tau^{2+\eta}) +
\mathcal{O}( \tau^3 h^{\alpha_1}) +  \mathcal{O}( \tau^4
h^{\beta_2})\right).
\end{array}
\end{array} $$
This provides the bound for the local errors $l_n^{[2]}$. The rest
of the proof follows as in the previous theorems.  \hfill $\Box$

%%%%%%%%%%%%%%%%%%%%%%%%%%%%%%%%%%%%%%%%%%%%%%
%\input{sect3-PRE}
\section{Application of the convergence results for Dirichlet Boundary Conditions in parabolic problems}

Let us next consider the $2D$ semi-linear diffusion-reaction model
($\varepsilon$ is a positive constant)
\begin{equation}\label{2D-dif-reac}
u_t=\varepsilon(u_{xx}+u_{yy}) + g(x,y,t),\; (x,y)\in
(0,1)^2,\,t\in[0,1],\;\varepsilon>0,
\end{equation}
with prescribed Dirichlet boundary conditions and an initial
condition. The PDE  is discretized on  uniform spatial meshes
$(x_i,y_j)=(ih,jh)$, $h= N^{-1}$, $1\le i,j\le N-1$, where $N-1$ is
the number of interior grid-points for each spatial variable. We
shall assume that the exact solution of the PDE (\ref{2D-dif-reac})
is regular enough when  $(x,y,t)\in [0,1]^2\times [0,t^*]$.
 Let us denote
$u_h(t):=(u_{i,j}(t))_{i,j=1}^{N-1}$ with a row-wise ordering, where
$u_{i,j}(t):=u(x_i,y_j,t)$ for $0\leq i,j\leq N$. Then, by using
second-order central differences, we obtain for the exact solution
of (\ref{2D-dif-reac}) on the grid a semi-discrete system
(\ref{pde}) with dimension $m=(N-1)^2$
\begin{equation}\label{semilin-exact-2D}
u_h'(t)=\varepsilon J u_h(t)+g_h(t)+\sigma_h(t)+\varepsilon
 h^{-2}u_{\Gamma_h}(t),
\end{equation}
where
\begin{equation}\label{J-2D}\begin{array}{c}
J:=J_1+J_2, \;\;J_1=I_{N-1}\otimes B_{N-1},\;\; J_2=B_{N-1}\otimes
I_{N-1},\\ B_{N-1}=h^{-2}TriDiag(1,-2,1)\in
\mathbb{R}^{(N-1)\times(N-1)},\quad h=1/N. \end{array}
\end{equation}
Moreover, $g_h(t)=(g(x_i,y_j,t))_{i,j=1}^{N-1}$, $\Vert \sigma_h(t)
\Vert_{2,h} =\mathcal{O}(h^2)$ ($0 \le t \le t^*$), whereas
$u_{\Gamma_h}(t)$ contains the values of the exact solution on the
boundary, i.e.,
\begin{equation}\label{uh-boundary}
u_{\Gamma_h}(t)=u_h^{(0,y)}(t)\otimes e_1+u_h^{(1,y)}(t)\otimes
e_{N-1}+e_1\otimes u_h^{(x,0)}(t)+e_{N-1}\otimes u_h^{(x,1)}(t),
\end{equation}
with $u_h^{(0,y)}(t)=(u_{0,j}(t))_{j=1}^{N-1}$,
$u_h^{(1,y)}(t)=(u_{N,j}(t))_{j=1}^{N-1}$,
$u_h^{(x,0)}(t)=(u_{i,0}(t))_{i=1}^{N-1}$ and
$u_h^{(x,1)}(t)=(u_{i,N}(t))_{i=1}^{N-1}$. Above,
$\{e_1,\ldots,e_{N-1}\}$ denotes the canonical
basis in $\mathbb{R}^{N-1}$.

For the proof of the convergence results we need the  lemma
\ref{lema-valP2} and  the lemma \ref{lema-Pgeneral} given below.
These lemmas can be derived from the material in \cite[pp.
96-300]{HV} (see from Lemma 6.1 to Lemma 6.5).  Lemmas
\ref{lema-valP2} and \ref{lema-Pgeneral} supply sharp values for the
constants $\alpha_l,\;\beta_l$ and $\eta$ appearing in the {\bf
P}-assumptions of section 3. These constants together with the
convergence theorems provide specific orders of convergence of the
MoL approach for several AMF$_q$-RK methods, in particular for the
{\sf AMF$_q$-Rad} methods presented in section 2.

The norm considered here for vectors, is the weighed Euclidean norm
$$ \displaystyle{\Vert (v_{ij})_{i,j=1}^{N-1}\Vert_{2,h}
:=\sqrt{\frac{1}{N^2}\sum_{i,j=1}^{N-1}|v_{ij}|^2}=h\Vert
(v_{ij})_{ij=1}^{N-1}\Vert_2, }$$ and for matrices the corresponding
operator norm.

\begin{lemma}\label{lema-valP2}
Assume that exact solution $u(x,y,t)$ of the 2D-PDE problem
(\ref{2D-dif-reac}) has as many continuous partial derivatives as
needed in the analysis   in $(x,y,t)\in [0,1]^2\times [0,t^*]$. Then
for  $ k=1,2\ldots$ and $\omega<\frac{1}{4}$ we have that,
$$\begin{array}{rcl}\norm{J^{\omega}u_h^{(k)}(t)}_{2,h}&=&\mathcal{O}(1),
\quad \mbox{\rm and moreover } \\
\norm{J^{1+\omega}u_h^{(k)}(t)}_{2,h}&=&\mathcal{O}(1), \; \mbox{\rm
whenever } u^{(1)}_{\Gamma_h}(t)\equiv 0.
\end{array}$$
\end{lemma}

\begin{lemma}\label{lema-Pgeneral}
Assume that exact solution $u(x,y,t)$ of the 2D-PDE problem
(\ref{2D-dif-reac}) has as many continuous partial derivatives as
needed in the analysis   in $(x,y,t)\in [0,1]^2\times [0,t^*]$.
Then, for $l=0,1,...$ we have that,
\begin{equation}\label{P-hypo-gen}
\norm{(P-J)^l u_h^{(k)}(t)}_{2,h}=\mathcal{O}(\tau^l
h^{\alpha_l}),\qquad \norm{J^{-1}(P-J)^l
u_h^{(k)}(t)}_{2,h}=\mathcal{O}(\tau^l h^{\beta_l}),
\end{equation}
where
\begin{equation}\label{P-hypo-alpha}
\alpha_l= \left\{\begin{array}{ll}-\max\{0,3+4(l-2)\}, &\quad {\rm
if} \; u^{(1)}_{\Gamma_h}(t)\equiv 0,\\-\max\{0,3+4(l-1)\}, & \quad
{\rm otherwise},\end{array}\right.
\end{equation}
and
\begin{equation}\label{P-hypo-beta}
\beta_l= \left\{\begin{array}{ll}-\max\{0,1+4(l-2)\}, &\quad {\rm
if} \; u^{(1)}_{\Gamma_h}(t)\equiv 0,\\-\max\{0,1+4(l-1)\}, & \quad
{\rm otherwise}.\end{array}\right.
\end{equation}
\end{lemma} \hfill $\Box$

We next give a convergence theorem for 2D-parabolic PDEs when the
MoL approach with   {\sf AMF$_q$-Rad} methods in section 2 are
applied to the time discretization. The results still hold for
3D-parabolic problems (even for $d$D-parabolic problems and $d\ge
3$) and Time-Independent Dirichlet boundary conditions, but the
proof requires some extra length to be included here.

\begin{theorem}\label{sev-th-7} The global errors (GE) in the weighted Euclidean norm  of the MoL
approach for the  2D-PDE (\ref{2D-dif-reac}) when the spatial
semi-discretization is carried out with second order central
differences and the  time integration is performed with AMF$_q$-RK
methods, are given in Table  \ref{estimates1D}. There,
$\varrho=\min\{1,\tau^2h^{-1}\}$ and $\mathcal{O}(\tau^{2.25^*})$ is
meant for $\mathcal{O}(\tau^{\mu})$ where  $\mu<2.25$ is any
constant.
\begin{table}[h!]
\centering
\begin{tabular}{|c|c|c|}  \hline
$(\tau\rightarrow 0^+,\:h\rightarrow 0^+)$ & GE (Time-Indep.) &  GE (Time-Dep.)  \\[0.2pc]
 \hline
{\sf AMF$_1$-Rad} & $\mathcal{O}(h^2)+\mathcal{O}(\tau^2)$ &
$\mathcal{O}(h^2) + \mathcal{O}(\varrho)$ \\[0.2pc]\hline {\sf AMF$_2$-Rad} &
$\mathcal{O}(h^2)+\mathcal{O}(\tau^3)+\tau^2\mathcal{O}(\varrho)$ &
$\mathcal{O}(h^2) + \mathcal{O}(\varrho)$ \\[0.2pc]\hline {\sf AMF$_3$-Rad} &
$\mathcal{O}(h^2)+\mathcal{O}(\tau^{2.25^*})$ &
$\mathcal{O}(h^2)+\mathcal{O}(\varrho)$ \\[0.2pc]
\hline
\end{tabular}\caption{\scriptsize Global error estimates in the weighted Euclidean norm
for  Time-Dependent Dirichlet boundary conditions (in short
Time-Dep.) and  Time-Independent Dirichlet boundary conditions (in
short Time-Indep.).} \label{estimates1D}
\end{table}

\end{theorem}

{\bf Proof.} In all cases we have that the stage order of the
underlying Runge-Kutta Radau IIA method is $\ell=2$ and  the order
of the  spatial semi-discretization is $r=2$. Moreover, all the
three methods {\sf AMF$_q$-Rad} ($q=1,\:2,\:3$) are
A($\pi/2$)-stable for a 2-splitting as it is shown in
\cite{sevedom-AMFestab} for the cases $q=1$ and $q=2$ and in
\cite{apnum-sevsole10} for the case $q=3$. Also, it should be
noticed that (\ref{spect-1}) holds.

We start with the {\sf AMF$_1$-Rad} method. We  have for the case of
Time-Independent Dirichlet Boundary conditions that the derivative
regarding $t$ vanishes on boundary points  $(x,y)\in \Gamma_h$, i.e.
$ u^{(1)}_{\Gamma_h}(t)\equiv 0$. From Lemma \ref{lema-Pgeneral} we
get that $\alpha_1=0$ and $\beta_1=0$. Then the bound for the global
errors follows from Theorem \ref{sev-th-3}. For the case of
Time-Dependent Dirichlet Boundary conditions, from Lemma
\ref{lema-Pgeneral}, we have that $\alpha_1=-3$ and $\beta_1=-1$.
Then,  the bound for the global errors follows from Theorem
\ref{sev-th-3}. The bound also applies to the {\sf AMF$_2$-Rad}
method for Time-Dependent Dirichlet BCs, because this method fulfils
the assumptions in  Theorem \ref{sev-th-3}.

For the case of the {\sf AMF$_2$-Rad} method  and Time-Independent
Dirichlet BCs we apply Theorem \ref{sev-th-3} for the case
$\varrho=1$ and  Theorem \ref{sev-th-5} with $l^*=1$ for the case
$\varrho=\tau^2h^{-1}$. Observe that from Lemma \ref{lema-Pgeneral}
we have that $\alpha_1=0$ and $\beta_1=0$ and $\beta_2=-1$. Moreover
the {\sf AMF$_2$-Rad} method fulfils all the assumptions in  Theorem
\ref{sev-th-5} by taking  $\eta=1$, see also Lemma \ref{lema-valP2}.

For the case of the {\sf AMF$_3$-Rad} method  and Time-Independent
Dirichlet BCs we apply Theorem \ref{sev-th-4} with any $\eta<1.25$,
see Lemma \ref{lema-valP2}. Observe that in this case
$\alpha_1=0,\;\beta_1=0$. Then from Theorem \ref{sev-th-4} the
global errors are of size $\mathcal{O}(h^2)+\mathcal{O}(\tau^{2})$.
The proof  that the order can be increased up to
$\mathcal{O}(h^2)+\mathcal{O}(\tau^{2.25^*})$ requires some extra
technical details that we have omitted for simplicity. The case of
Time-Dependent Dirichlet BCs follows from Theorem \ref{sev-th-4}
too, but in this case $\beta_1=-1$. \hfill $\Box$

%%%%%%%%%%%%%%%%%%%%%%%%%%%%%%%%%%%%%%%%%%%%%%%%%%
\subsection{Numerical Experiments}
We have performed some numerical experiments on two 2D-PDE and
3D-PDE problems of parabolic type in order to illustrate the
convergence results presented in former sections for the {\sf
AMF$_q$-Rad} methods.

\begin{enumerate}
\item {\sf Problem 1} is the 2D-PDE problem (\ref{2D-dif-reac}) with
diffusion parameter $\varepsilon=0.1$ and Dirichlet Boundary
Conditions and an Initial Condition so that
\begin{equation}\label{2D-solution} u(x,y,t)= 10x(1-x)y(1-y)e^t + \beta
e^{2x-y-t}, \end{equation} is the exact solution. The case $\beta=0$
provides Time-Independent Boundary conditions and no spatial error
($\sigma_h(t)\equiv 0$, due to the polynomial nature of the exact
solution). The case $\beta=1$ provides Time-Dependent boundary
conditions and spatial discretizations errors of order two.
\item {\sf Problem 2} is the 3D-PDE problem (\ref{3D-dif-reac}) with
diffusion parameter $\varepsilon=0.1$
\begin{equation}\label{3D-dif-reac}\begin{array}{c}
u_t(\overrightarrow{x},t)=  \varepsilon\:
\Delta u(\overrightarrow{x},t)+g(\overrightarrow{x},t),\\
t\in[0,1],\;\;\overrightarrow{x}=(x,y,z)\in (0,1)^3  \in
\mathbb{R}^3,
\end{array}
 \end{equation}
and Dirichlet Boundary Conditions and an Initial Condition so that
\begin{equation}\label{2D-solution} u(x,y,t)= 64x(1-x)y(1-y)z(1-z)e^t + \beta
e^{2x-y-z-t}, \end{equation} is the exact solution. Again, the case
$\beta=0$ provides Time-Independent Boundary conditions and no
spatial error  and the case $\beta\ne 0$ provides Time-Dependent
boundary conditions and spatial discretizations errors of order two.
\end{enumerate}
%%%%%%%%%%%%%
On the end-point of the time interval $t^*=1$,   in the weighted
Euclidean norm we have computed as specified in (\ref{sev-equa-1}),
the global errors $\epsilon_2(h,\tau)$ ($y_{\rm met}(t^*)$ denotes
the numerical solution at $t^*$ by the method considered), the
number of significant figures of the global errors
$\delta_2(h,\tau)$ and the estimated order of the global errors
$p(h,\tau)$  as powers of $h$ when $r=\tau/h$ is kept constant and
both $\tau$ and $h$ tend to zero.
\begin{equation}\label{sev-equa-1}\begin{array}{c}
\epsilon_2(h,\tau):=\norm{u_h(t^*)-y_{\rm met}(t^*)}_{2,h}, \quad
\delta_2(h,\tau)=-\log_{10} \epsilon_2(h,\tau)\\
p(h,\tau)= (\delta_2(h/2,\tau/2)- \delta_2(h,\tau))/\log_{10}
2.
\end{array}
\end{equation}

In the Tables \ref{table-linear2D-1}, \ref{table-linear2D-2} and
\ref{table-linear3D-1}  we have considered  for each $h$ the
time-stepsize $\tau=q h$ for the corresponding {\sf AMF$_q$-Rad}
method ($q=1,2,3$), so that all the methods make use of the same
number of $f$-evaluations and similar CPU times in the computations.
In those tables we have displayed the number of significant figures
in the global errors  $\delta_2(h,\tau)$ and in brackets the
estimated orders $p(h,\tau)$ of each method.

From Theorem \ref{sev-th-7},  the global errors are expected to be
of size $h^\mu$ (observe that $\tau/h$ is kept  constant) where:
\begin{enumerate}
\item for the {\sf AMF$_1$-Rad} method, $\mu=2$ if  Time-Independent BCs are considered and $\mu=1$ if  Time-Dependent
BCs are imposed. This nicely fits with the results displayed in
Table \ref{table-linear2D-1} (Time-Independent BCs) and in Table
\ref{table-linear2D-2} (Time-Dependent BCs) for the 2D-PDE problem.
Moreover, the convergence order is still $\mu=2$ in the 3D-PDE
problem for Time-Independent BCs as it can be seen in Table
\ref{table-linear3D-1}.
\item For the {\sf AMF$_2$-Rad} method, $\mu=3$ if  Time-Independent BCs are considered and $\mu=1$ if  Time-Dependent
BCs are imposed. This  fits well with the results displayed in Table
\ref{table-linear2D-1} (Time-Independent BCs) and in Table
\ref{table-linear2D-2} (Time-Dependent BCs) for the 2D-PDE problem.
Moreover, the convergence order is also $\mu=3$ in the 3D-PDE
problem for Time-Independent BCs as it can be observed  in Table
\ref{table-linear3D-1}.
\item For the {\sf AMF$_3$-Rad} method, $\mu=2.25^*$ if  Time-Independent BCs are considered
and $\mu=1$ if  Time-Dependent BCs are imposed. This  can be
observed in Table \ref{table-linear2D-1} (Time-Independent BCs) and
in Table \ref{table-linear2D-2} (Time-Dependent BCs) for the 2D-PDE
problem. Moreover, the convergence order also approaches  to
$\mu=2.3$ in the 3D-PDE problem for Time-Independent BCs as shown in
Table \ref{table-linear3D-1}.
\end{enumerate}

%%%%%%%%%%%%%%%%%%%%%%%%%%%%%%
%Tables
\begin{table}[h!]
\centering
\begin{tabular}{|l|l|l|l|}  \hline
 $h$ & $\begin{array}{c}\mbox{\sf AMF$_1$-Rad} \;({p})\\ \tau/h= 1\end{array}$ &
 $\begin{array}{c}\mbox{\sf  AMF$_2$-Rad} \;({p})\\ \tau/h= 2\end{array}$ &
  $\begin{array}{c}\mbox{\sf  AMF$_3$-Rad} \;({p})\\ \tau/h= 3\end{array}$   \\ \hline
$\frac{1}{24}$ & $ \delta_2=3.74 \; (2.03)$ & $ \delta_2=4.94 \;
(2.82)$ & $ \delta_2=4.90 \; (3.56)$
\\\hline
$\frac{1}{48}$ & $ \delta_2=4.35 \; (2.03)$ & $ \delta_2=5.79 \;
(2.89)$ & $ \delta_2=5.67 \; (2.42)$
\\\hline$\frac{1}{96}$ & $ \delta_2=4.96
\; (1.99)$ & $ \delta_2=6.66 \; (2.92)$ & $ \delta_2=6.40 \; (2.36)$
\\\hline$\frac{1}{192}$ & $ \delta_2=5.56
\; (1.99)$ & $ \delta_2=7.54 \; (2.93)$ & $ \delta_2=7.11 \; (2.29)$
\\\hline$\frac{1}{384}$ & $ \delta_2=6.16
\; (2.03)$ & $ \delta_2=8.42 \; (2.96)$ & $ \delta_2=7.80 \; (2.29)$
\\\hline$\frac{1}{768}$ & $ \delta_2=6.77
\; (--)$ & $ \delta_2=9.31 \; (--)$ & $ \delta_2=8.49 \; (--)$
\\\hline
\end{tabular}
\caption{\scriptsize Significant correct digits  ($l_{2,h}$-norm)
for the 2D-PDE problem with Time-Independent Dirichlet BCs
($\beta=0$). In brackets the estimated orders of convergence (by
halving both the spatial resolution $h$ and the time-stepizes $\tau$
and taking ratio $r=\tau/h$).}\label{table-linear2D-1}
\end{table}

\begin{table}[h!]
\centering
\begin{tabular}{|l|l|l|l|}  \hline
 $h$ & $\begin{array}{c}\mbox{\sf AMF$_1$-Rad} \;({p})\\ \tau/h= 1\end{array}$ &
 $\begin{array}{c}\mbox{\sf  AMF$_2$-Rad} \;({p})\\ \tau/h= 2\end{array}$ &
  $\begin{array}{c}\mbox{\sf  AMF$_3$-Rad} \;({p})\\ \tau/h= 3\end{array}$   \\ \hline
$\frac{1}{24}$ & $ \delta_2=3.02 \; (1.00)$ & $ \delta_2=2.79 \;
(0.76)$ & $ \delta_2=2.52 \; (0.66)$
\\\hline
$\frac{1}{48}$ & $ \delta_2=3.32 \; (0.97)$ & $ \delta_2=3.02 \;
(0.83)$ & $ \delta_2=2.72 \; (0.76)$
\\\hline$\frac{1}{96}$ & $ \delta_2=3.61
\; (1.00)$ & $ \delta_2=3.27 \; (0.90)$ & $ \delta_2=2.95 \; (0.86)$
\\\hline$\frac{1}{192}$ & $ \delta_2=3.91
\; (1.00)$ & $ \delta_2=3.54 \; (0.93)$ & $ \delta_2=3.21 \; (0.91)$
\\\hline$\frac{1}{384}$ & $ \delta_2=4.21
\; (1.03)$ & $ \delta_2=3.82 \; (0.97)$ & $ \delta_2=3.48 \; (0.97)$
\\\hline$\frac{1}{768}$ & $ \delta_2=4.52
\; (--)$ & $ \delta_2=4.11 \; (--)$ & $ \delta_2=3.77 \; (--)$
\\\hline
\end{tabular}
\caption{\scriptsize Significant correct digits  ($l_{2,h}$-norm)
for the 2D-PDE problem with Time-Dependent Dirichlet BCs
($\beta=1$). In brackets the estimated orders of convergence (by
halving both the spatial resolution $h$ and the time-stepizes $\tau$
and taking ratio $r=\tau/h$).}\label{table-linear2D-2}
\end{table}

\begin{table}[h!]
\centering
\begin{tabular}{|l|l|l|l|}  \hline
 $h$ & $\begin{array}{c}\mbox{\sf AMF$_1$-Rad} \;({p})\\ \tau/h= 1\end{array}$ &
 $\begin{array}{c}\mbox{\sf  AMF$_2$-Rad} \;({p})\\ \tau/h= 2\end{array}$ &
  $\begin{array}{c}\mbox{\sf  AMF$_3$-Rad} \;({p})\\ \tau/h= 3\end{array}$   \\ \hline
$\frac{1}{24}$ & $ \delta_2=3.40 \; (2.03)$ & $ \delta_2=4.31 \;
(2.96)$ & $ \delta_2=4.53 \; (2.69)$
\\\hline
$\frac{1}{48}$ & $ \delta_2=4.01 \; (2.03)$ & $ \delta_2=5.20 \;
(2.96)$ & $ \delta_2=5.34 \; (2.59)$
\\\hline$\frac{1}{96}$ & $ \delta_2=4.62
\; (--)$ & $ \delta_2=6.09 \; (--)$ & $ \delta_2=6.12 \; (--)$
\\\hline
\end{tabular}
\caption{\scriptsize Significant correct digits  ($l_{2,h}$-norm)
for the 3D-PDE problem with Time-Independent Dirichlet BCs
($\beta=0$). In brackets the estimated orders of convergence (by
halving both the spatial resolution $h$ and the time-stepizes $\tau$
and taking ratio $r=\tau/h$).}\label{table-linear3D-1}
\end{table}

As a conclusion we can say that the convergence results presented in
Theorem \ref{sev-th-7}  seem to be sharp for 2D-parabolic problems
and that they still hold for $d$D-parabolic problems ($d > 2$) when
Time-Independent boundary conditions are considered. The proof of
this fact requires some additional work and is not presented here.
On the other hand, the convergence results are very poor when
Time-Dependent Boundary conditions are considered. However, in such
a situation we have  developed a very simple technique (Boundary
Correction Technique) to recover the convergence order as if
Time-Independent Boundary conditions were considered. The
explanation of the Boundary Correction Technique and the proof of
the convergence orders  requires some extra length and will be the
objective of another paper.

It is also important to remark that although we have considered in
Theorem \ref{sev-th-7}, second-order central differences for the
spatial discretization, the convergence results also hold for most
of the usual spatial discretizations as long as they are stable and
consistent with order $r\ge 1$. Numerical experiments carried by the authors 
seem to indicate that the convergence results also hold for many classes of non-linear problems.

%%%%%%%%%%%%%%%%%%%%%%%%%%%%%%%%%%%%%%%%%%%%%%%%%%%%%%%%%%%%%%%%%%%%%
%\input{references-PRE}

\end{document}